
\input psfig
\input amssym.def
\input amssym
\magnification=1100
\baselineskip = 0.25truein
\lineskiplimit = 0.01truein
\lineskip = 0.01truein
\vsize = 8.5truein
\voffset = 0.2truein
\parskip = 0.10truein
\parindent = 0.3truein
\settabs 12 \columns
\hsize = 5.4truein
\hoffset = 0.4truein

\setbox\strutbox=\hbox{%
\vrule height .708\baselineskip
depth .292\baselineskip
width 0pt}
\font\caps=cmcsc10

\font\bigtenrm=cmr10 at 14pt
\font\bigmaths=cmmi10 at 16pt

\def\sqr#1#2{{\vcenter{\vbox{\hrule height.#2pt
\hbox{\vrule width.#2pt height#1pt \kern#1pt
\vrule width.#2pt}
\hrule height.#2pt}}}}
\def\square{\mathchoice\sqr46\sqr46\sqr{3.1}6\sqr{2.3}4}

\centerline{\bigtenrm HEEGAARD SPLITTINGS,}
\centerline{\bigtenrm THE VIRTUALLY HAKEN CONJECTURE}
\centerline{\bigtenrm \textfont1=\bigmaths
AND PROPERTY ($\tau$)}
\footnote{}{MSC (2000): 57N10, 57M10 (primary);
58C40, 05C25, 57M15 (secondary)}

\tenrm
\vskip 14pt
\centerline{MARC LACKENBY}
\vskip 18pt

\tenrm
\centerline{\caps 1. Introduction}
\vskip 6pt

The behaviour of finite covers of 3-manifolds is a
very important, but poorly understood, topic.
There are three, increasingly strong, conjectures in the field
that have remained open for over twenty years
-- the virtually Haken conjecture, the positive virtual $b_1$ conjecture
and the virtually fibred conjecture. Any of these would have
profound ramifications for 3-manifold theory. In this paper,
we explore the interaction of these conjectures with
the following seemingly unrelated areas: eigenvalues of
the Laplacian, and Heegaard splittings. 

We first give a 
necessary and sufficient condition, in terms of spectral
geometry, for a finitely presented group to have a finite
index subgroup with infinite abelianisation. 
This result in geometric group theory may have
uses beyond 3-manifold theory.
We also show that, for negatively curved 3-manifolds,
this is equivalent to a statement about generalised
Heegaard splittings. This provides one link between
the positive virtual $b_1$ conjecture, Heegaard splittings
and the Laplacian.

In a second direction, we define a new invariant of
3-manifolds: their Heegaard gradient. This measures
the growth rate of the Heegaard genus of the manifold's finite
covering spaces, as a function of their degree.
We formulate a conjecture about Heegaard gradient,
and provide some supporting evidence. We show
that this, together with a conjecture
of Lubotzky and Sarnak about Property ($\tau$), would
imply the virtually Haken conjecture for hyperbolic 3-manifolds.
Property ($\tau$) is a concept due to Lubotzky and
Zimmer [37], which encodes the behaviour of the first eigenvalue
of the Laplacian in finite covering spaces.

Along the way, we prove a number of unexpected theorems
about 3-manifolds. For example, we show that
for any closed 3-manifold that fibres over the circle with
pseudo-Anosov monodromy, any cyclic cover dual to the fibre 
of sufficiently large degree has an irreducible, weakly reducible
Heegaard splitting. Also, we establish lower and upper bounds on 
the Heegaard genus of the congruence covers of an arithmetic hyperbolic
3-manifold, which are linear in volume. 

\vskip 6pt
\noindent {\caps The virtually Haken conjecture and related conjectures}
\vskip 6pt

The following will motivate much of this paper.

\noindent {\bf Virtually Haken conjecture.} {\sl A compact orientable
irreducible 3-manifold with infinite fundamental group is
virtually Haken.}

This is important for many reasons. Haken manifolds have played a
central r\^ole in 3-manifold theory for the past 40 years. Although
many closed orientable irreducible 3-manifolds with infinite fundamental group are 
known to be non-Haken, it would be useful and satisfying to know that
any such 3-manifold $M$ was finitely covered by a Haken 3-manifold
$\tilde M$. This would have many implications for $M$. Firstly,
the embedded incompressible surface in $\tilde M$ would project
to a $\pi_1$-injective immersed surface in $M$. Thus, $\pi_1(M)$ would contain the
fundamental group of a closed orientable surface with positive
genus. Secondly, many of the important properties that are known
to hold for $\tilde M$ descend to properties of $M$. A notable
example of this phenomenon is the theorem that virtually
Haken 3-manifolds satisfy the geometrisation conjecture. For,
Thurston showed that if a 3-manifold is Haken, then it
satisfies the geometrisation conjecture [40] and
it follows from the equivariant sphere theorem [38], 
the Seifert fibre space theorem ([7],[18],[54]) and the fact
that virtually hyperbolic 3-manifolds are hyperbolic ([19],[20]) that
if a compact orientable irreducible 3-manifold with infinite fundamental group is
finitely covered by a 3-manifold satisfying the geometrisation
conjecture, then it does also. There is the possibility that the geometrisation
conjecture has recently been proved by other means, due to the work
of Hamilton ([22], [23], [24], [25]) and 
Perelman ([42], [43], [44]). This is fortuitous, since
there seems to be little hope of proving the virtually Haken
conjecture without assuming that the manifold is geometric, or
making some similar hypothesis. A solution to the virtually Haken conjecture
just for hyperbolic 3-manifolds would be a considerable achievement
in its own right.

The following implies the virtually Haken conjecture.

\noindent {\bf Positive virtual $b_1$ conjecture.} {\sl 
A compact orientable irreducible atoroidal 3-manifold with infinite 
fundamental group has a finite cover with positive first
Betti number.}

A recent analysis by Dunfield and Thurston
of the 10986 closed hyperbolic 3-manifolds in the Hodgson-Weeks census
showed that they all satisfy this conjecture [13].

The following conjecture, due to Thurston [58], is stronger still.

\noindent {\bf Virtually fibred conjecture.} {\sl 
A compact orientable irreducible 3-manifold with boundary 
a (possibly empty) collection of tori, whose fundamental 
group is infinite and contains no ${\Bbb Z} \times {\Bbb Z}$
subgroup, is finitely covered by a surface bundle over the circle.}

This seems less likely to be true, since there are very few 
known examples of 3-manifolds that are virtually fibred but not
fibred ([48], [33]).

\vskip 6pt
\noindent{\caps Heegaard splittings}
\vskip 6pt

Casson and Gordon [6] introduced the notion of a weakly reducible
Heegaard splitting and showed that if a compact orientable 3-manifold admits
a splitting that is irreducible but weakly reducible, then it
is Haken. This technology was developed by Scharlemann
and Thompson [51], who introduced generalised Heegaard splittings,
which are a decomposition of the 3-manifold along closed
surfaces, together with Heegaard surfaces for the complementary
regions. They defined a complexity for such splittings, 
and termed a generalised Heegaard splitting
{\sl thin} when it has minimal complexity. They showed that
if an irreducible 3-manifold has a thin generalised Heegaard splitting that is 
not a Heegaard splitting, then the manifold contains a
closed essential surface. In particular, it is Haken.

Our approach in this paper is to study the behaviour of Heegaard genus
under finite covers. In fact, it is slightly more convenient
to use the Euler characteristic of Heegaard surfaces. We therefore
define the {\sl Heegaard Euler characteristic} $\chi_-^h(M)$
of a compact orientable 3-manifold $M$ to be the negative
of the maximal Euler characteristic of a Heegaard surface. 
Of course, since Heegaard surfaces are always connected,
$\chi_-^h(M)$ is linearly related to the more familiar 
Heegaard genus $g(M)$ by the formula $\chi_-^h(M) = 2g(M) - 2$.
We also define the {\sl strong Heegaard Euler characteristic}
$\chi_-^{sh}(M)$ of $M$
to be the negative of the maximal Euler characteristic
of a strongly irreducible Heegaard surface, or infinity
if $M$ does not have such a Heegaard surface.
One can study many interesting asymptotic quantities associated
with these invariants, but we focus on two.

\noindent {\bf Definition.} Let $\{ M_i \rightarrow M \}$ be a
collection of
finite covers of a compact orientable 3-manifold $M$, with
degree $d_i$. The {\sl infimal Heegaard gradient} of $\{ M_i \rightarrow
M \}$
is 
$$\inf_i {\chi_-^h(M_i) \over d_i }.$$
The {\sl infimal strong Heegaard gradient} of
$\{ M_i \rightarrow M \}$ is
$$\liminf_i {\chi_-^{sh}(M_i) \over d_i },$$
which may be infinite.
When the covers $\{ M_i \rightarrow M \}$ are not
mentioned, they are understood to be all the finite
covers of $M$. For brevity, we sometimes drop the word `infimal'.

A Heegaard splitting for $M$ lifts to one for $M_i$, and its
Euler characteristic is scaled by $d_i$. Thus, the Heegaard gradient
of $M$ is at most $\chi_-^h(M)$. Heegaard gradient therefore
measures the degeneration of Heegaard Euler characteristic
in finite sheeted covers. There are variants of the
above definitions, where $\inf$ is replaced by $\liminf$,
and where only regular covers are considered. We explore these
in \S3.

There exist hyperbolic 3-manifolds with zero Heegaard gradient.
For example, $\chi_-^h(M)$ of a closed orientable 3-manifold $M$ that fibres 
over the circle with fibre $F$ is at most $2|\chi(F)|+4$.
Thus, by passing to cyclic covers $M_i$, we can
make $d_i$ arbitrarily large, but keep $\chi_-^h(M_i)$ bounded above.
So, if a 3-manifold virtually fibres over
the circle, then its infimal Heegaard gradient is zero. The following 
conjecture proposes that this is the only source of this phenomenon
for hyperbolic 3-manifolds.

\noindent {\bf Heegaard gradient conjecture.} {\sl A compact orientable 
hyperbolic 3-manifold has zero Heegaard gradient if and only if it
virtually fibres over the circle.}

Along similar lines, we also propose the following.

\noindent {\bf Strong Heegaard gradient conjecture.} 
{\sl Any closed orientable hyperbolic
3-manifold has positive strong Heegaard gradient.}

As we shall see, establishing either of these conjectures
would be a step towards proving the virtually Haken conjecture
for hyperbolic 3-manifolds.
We will provide some supporting evidence for them in this paper.

\vskip 6pt
\noindent {\caps Property ($\tau$)}
\vskip 6pt

Property ($\tau$) for finitely generated groups was 
defined by Lubotzky and Zimmer [37]. It has a number of equivalent
definitions, relating to graph theory, differential geometry,
and group representations, which are recalled in \S2.
In this paper, we will focus on its interpretation
in terms of the smallest eigenvalue of the Laplacian,
and the Cheeger constant.
The following conjecture is central to this paper.

\noindent {\bf Conjecture.} [Lubotzky-Sarnak] {\sl The fundamental
group of a finite volume hyperbolic 3-manifold fails 
to have Property ($\tau$).}

Evidence for this conjecture is presented in [36]. Its
motivation arises from the similarity between Property 
($\tau$) and Kazhdan's Property (T) [30], for which the
corresponding assertion is true. Many finitely generated
groups are known not to have Property ($\tau$). For example,
a finitely generated group with a finite index subgroup
having infinite abelianisation does not have Property ($\tau$).
Hence, the positive virtual $b_1$ conjecture implies
the Lubotzky-Sarnak conjecture. Another class of groups
that fail to have Property ($\tau$) are infinite, amenable, 
residually finite groups.

\vskip 6pt
\noindent {\caps The main results}
\vskip 6pt

In \S2, we establish a necessary and sufficient condition for
a finitely presented group to have a finite index subgroup
with infinite abelianisation, in terms
of eigenvalues of the Laplacian, and the geometry
of Schreier coset graphs. It is somewhat surprising that such
a characterisation should exist. 

\vfill\eject
\noindent {\bf Theorem 1.1.} {\sl Let $G$ be a finitely
presented group, and let $S$ be a finite set of generators. Let
$\{ G_i \}$ be its finite index subgroups, and
let $N(G_i)$ be the normaliser of $G_i$ in $G$. Let $X_i$
be the Schreier coset graph of $G/G_i$ induced by $S$. Then the
following are equivalent, and are independent of
the choice of $S$:
\item{1.} some $G_i$ has infinite abelianisation;
\item{2.} $G_i$ has infinite abelianisation for infinitely many $i$;
\item{3.} $\lambda_1(X_i) [G : N(G_i)]^2 [G:G_i]^2$
has a bounded subsequence;
\item{4.} $\lambda_1(X_i) [G : N(G_i)]^4 [N(G_i):G_i]$
has zero infimum;
\item{5.} $h(X_i) [G : N(G_i)] [G:G_i]$ has a bounded subsequence;
\item{6.} $h(X_i) [G : N(G_i)]^2 [N(G_i):G_i]^{1/2}$
has zero infimum.

\noindent Furthermore, if $G$ is the fundamental group
of some closed orientable Riemannian manifold $M$, and $M_i$ is
the cover of $M$ corresponding to $G_i$, then the above
are also equivalent to each of the following:
\item{7.} $\lambda_1(M_i) [G : N(G_i)]^2 [G:G_i]^2$
has a bounded subsequence;
\item{8.} $\lambda_1(M_i) [G : N(G_i)]^4 [N(G_i):G_i]$
has zero infimum;
\item{9.} $h(M_i) [G : N(G_i)] [G:G_i]$ has a bounded subsequence;
\item{10.} $h(M_i) [G : N(G_i)]^2 [N(G_i):G_i]^{1/2}$
has zero infimum.

}

For a finite graph or Riemannian manifold, $\lambda_1$
denotes the smallest non-zero eigenvalue of the Laplacian,
and $h$ is its Cheeger constant. The definitions of
these terms are recalled in \S2. 

Note that $[G:G_i]$ and $[N(G_i):G_i]$ have the
following simple topological interpretations: $[G:G_i]$
is the number of vertices of $X_i$; and $[N(G_i):G_i]$
is the number of covering translations of $X_i$.
Note also that if $G_i$ is a normal subgroup of $G$,
then $[G: N(G_i)]$ is one.

Theorem 1.1 should be
compared with the definition of Property $(\tau)$,
which is also recalled in \S2.

A sample application of the methods behind Theorem 1.1
gives the following result. Recall that a group presentation
is {\sl triangular} if each relation has length three.
Any finitely presented group has a finite triangular
presentation.

\noindent {\bf Theorem 1.2.} {\sl Let $X$ be the Cayley
graph of a finite group, arising from a finite triangular
group presentation. Then $h(X) \geq \sqrt{2/(3|V(X)|)}$.}

In \S3, we study generalised Heegaard splittings. In [51],
Scharlemann and Thompson defined a complexity
for a compact orientable 3-manifold $M$, which measures the minimal
complexity of a generalised Heegaard splitting.
It is a finite multi-set of integers. The largest
integer in this multi-set, minus 1, is denoted by $c_+(M)$. 
When $M$ is closed, irreducible, non-Haken and not $S^3$, this is
equal to $\chi_-^h(M)$ and $\chi_-^{sh}(M)$.
In all cases, other than when $M$ is $S^3$ or a 3-ball, $c_+(M) \leq \chi_-^h(M)$.
More details can be found in \S3.  In \S4,
we find an upper bound on the Cheeger constant for
a negatively curved 3-manifold $M$ in terms of $c_+(M)$.
Using Theorem 1.1,
this allows us reformulate the positive virtual $b_1$ conjecture
for negatively curved 3-manifolds in terms of generalised Heegaard splittings.

\noindent {\bf Theorem 1.3.} {\sl Let $M$ be a closed orientable
3-manifold with a negatively curved Riemannian metric. Let
$\{ M_i \rightarrow M \}$ be the finite covers of $M$. 
Then the following are equivalent,
and are each equivalent to the conditions in Theorem 1.1:
\item{1.} $b_1(M_i) > 0$ for infinitely many $i$;
\item{2.} $c_+(M_i) [\pi_1 M : N(\pi_1 M_i )] $ has a bounded subsequence;
\item{3.} $c_+(M_i) [\pi_1 M : N(\pi_1 M_i )] [N(\pi_1 M_i ):
\pi_1 M_i ]^{-1/2}$ has zero infimum.

}

Using (3) $\Rightarrow$ (1), this has the following immediate corollary.

\noindent {\bf Corollary 1.4.} {\sl Let $M$ be a 
closed orientable 3-manifold with a negatively curved
Riemannian metric. Suppose that, for some collection
of finite regular covers $\{ M_i \rightarrow M \}$ with
degree $d_i$, $\inf \chi_-^h(M_i) / \sqrt d_i = 0$.
Then $M$ satisfies the positive virtual
$b_1$ conjecture.}

As another consequence of the relation of the Cheeger
constant to Heegaard Euler characteristic, we obtain
the following.

\noindent {\bf Theorem 1.5.} {\sl Let $M$ be an
orientable 3-manifold that admits a complete, negatively
curved, finite volume Riemannian metric. Let $\{ M_i \rightarrow M \}$ be a
collection of finite covers of $M$. If
$\pi_1(M)$ has Property ($\tau$) with respect to
$\{ \pi_1(M_i) \}$, then the infimal Heegaard gradient of
$\{  M_i \rightarrow M \}$ is non-zero.}

The fundamental group of an arithmetic hyperbolic
3-manifold has Property ($\tau$) with respect to its
congruence subgroups [35]. Hence, we have the following
immediate corollary.

\noindent {\bf Corollary 1.6.} {\sl Let $M$ an arithmetic
hyperbolic 3-manifold. Then there are positive constants $c$ and $C$,
such that for any congruence cover $M_i \rightarrow M$,
$$c \ {\rm Volume}(M_i) \leq \chi_-^h(M_i) \leq C
\ {\rm Volume}(M_i).$$

}

This may be viewed as supporting evidence for
the strong Heegaard gradient conjecture, since
it provides a large collection of finite covers with
non-zero Heegaard gradient and hence non-zero
strong Heegaard gradient.

In \S5, we relate Property ($\tau$), Heegaard splittings and
the virtually Haken conjecture, by proving the following
theorem. This may represent a viable approach to the
virtually Haken conjecture for hyperbolic 3-manifolds.

We term 
a collection of covers $\{ M_i \rightarrow M \}$ a {\sl lattice}
if, whenever $M_i \rightarrow M$ and $M_j \rightarrow M$ 
lies in the collection, so does
the cover corresponding to $\pi_1(M_i) \cap \pi_1(M_j)$.
For example, the set of all finite regular covers of $M$ forms a
lattice.

\noindent {\bf Theorem 1.7.} {\sl Let $M$ be a compact
orientable irreducible 3-manifold, with boundary a 
(possibly empty) collection of
tori. Let $\{ M_i \rightarrow M \}$ be a lattice
of finite regular covers of $M$. Suppose that
\item{1.} $\pi_1(M)$ fails to have Property ($\tau$) with
respect to $\{ \pi_1(M_i) \}$, and
\item{2.} $\{ M_i \rightarrow M \}$ has non-zero infimal
strong Heegaard gradient.

\noindent Then, for infinitely many $i$, $M_i$ has a
thin generalised Heegaard splitting which is not a Heegaard
splitting, and so contains a closed essential surface.
In particular, $M$ is virtually Haken.}

As a consequence, the virtually Haken conjecture for
finite volume hyperbolic 3-manifolds would follow
from the Lubotzky-Sarnak conjecture, together with
either the Heegaard gradient conjecture or the strong
Heegaard gradient conjecture.

An immediate consequence of Theorems 1.5 and 1.7 is the
following result.

\vfill\eject
\noindent {\bf Theorem 1.8.} {\sl Let $\{ M_i \rightarrow M \}$
be a lattice of finite regular non-Haken covers of a closed orientable negatively
curved 3-manifold $M$. Then the following are equivalent:
\item{1.} $\{M_i \rightarrow M \}$ has non-zero Heegaard gradient;
\item{2.} $\{M_i \rightarrow M \}$ has non-zero strong Heegaard gradient;
\item{3.} $\pi_1(M)$ has Property $(\tau)$ with respect to
$\{ \pi_1(M_i) \}$.}

It is striking that (3) above is a group-theoretic
property of $\pi_1(M)$ and $\{ \pi_1(M_i) \}$, 
whereas this is not at all obvious of either (1) or (2).
 
In \S6, we investigate the diameter of minimal surfaces in
negatively curved 3-manifolds. An immediate application
is the following, which provides some evidence for
the strong Heegaard gradient conjecture.

\noindent {\bf Theorem 1.9.} {\sl Let $M$ be a closed orientable
hyperbolic 3-manifold, and let $\{ M_i \rightarrow M \}$
be the cyclic covers dual to some non-trivial element
of $H_2(M)$. Then the strong Heegaard gradient of
$\{ M_i \rightarrow M \}$ is non-zero.}

Thus, in the above case, the covers $\{ M_i \rightarrow M \}$
satisfy each of the requirements of Theorem 1.7. Of course,
the final conclusion of Theorem 1.7 also holds trivially, but this
shows that the hypotheses of Theorem 1.7 are often satisfied.

Theorem 1.9 has the following interesting 
implication for Heegaard splittings of 
hyperbolic manifolds that fibre over the circle.
Rubinstein has independently proved a stronger version of
this result [50].

\noindent {\bf Corollary 1.10.} {\sl Let $M$ be a closed orientable
3-manifold that fibres over the circle with pseudo-Anosov
monodromy. Let $\{ M_i \rightarrow M \}$ be the cyclic covers dual to
the fibre. Then, for all but finitely many $i$, 
$M_i$ has an irreducible, weakly reducible, minimal genus 
Heegaard splitting.}

In \S7, we continue with the study of cyclic coverings.
The technical result in \S6 is used to prove the
following theorem, which supports the 
Heegaard gradient conjecture.

\noindent {\bf Theorem 1.11.} {\sl 
Let $M$ be a compact orientable finite volume hyperbolic 3-manifold, and
let $\{ M_i \rightarrow M \}$ be the cyclic covers dual to some
non-trivial element $z$ of $H_2(M, \partial M)$.
Then, the infimal Heegaard gradient of $\{ M_i \rightarrow M \}$
is zero if and only if $z$ is represented by a fibre.}

\vskip 6pt
\noindent {\caps Acknowledgements}
\vskip 6pt

I am indebted to Jason Manning, who pointed out an error
in an earlier version of Theorem 1.1. 
I thank Martin Lustig for a useful conversation. I thank 
Kazuhiro Ichihara whose work on the 
Heegaard gradient of Seifert fibre spaces [28] led me to
revise my initial version of the Heegaard gradient
conjecture. I would also like to thank the referee for
carefully reading this paper and for providing some
useful suggestions for improving it.

\vskip 18pt
\centerline{\caps 2. Positive virtual $b_1$ and spectral geometry} 
\vskip 6pt

In this section, we prove Theorem 1.1, which gives a number of equivalent
geometric characterisations of positive virtual $b_1$ for a
finitely presented group. It is stated in terms of the Cheeger
constant and the first eigenvalue of the Laplacian of
Cayley graphs and Riemannian manifolds. We now recall
the definitions of these terms.

Given a group $G$ with a finite set $S$ of generators,
and a subgroup $H$, the {\sl Schreier coset graph}
(or {\sl coset diagram}) of $G/H$ with respect to $S$ is
defined as follows. It has a vertex for each right coset
$Hg$, and vertices $Hg_1$ and $Hg_2$ are joined by
an oriented edge if and only if $g_2 = g_1 s$ for some
$s \in S$. When $H = 1$, this gives the {\sl Cayley graph} of
$G$ with respect to $S$. The Cayley graph of the cyclic group of
order 6, with respect to the generators $\{ 1,2 \}$,
is shown in Figure 1.

\vskip 18pt
\centerline{\psfig{figure=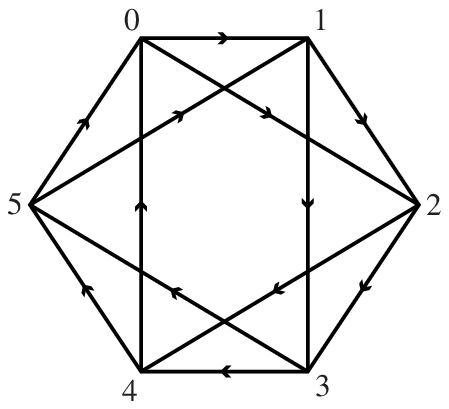,width=1.5in}}
\vskip 18pt
\centerline{Figure 1.}

Note that left multiplication by any element of $G$ induces an
automorphism of a Cayley graph for $G$.

Given a finite graph $X$, we denote its vertex set by $V(X)$
and its edge set by $E(X)$.
The {\sl Cheeger constant} of $X$, denoted $h(X)$,
is defined to be 
$$h(X) = \inf_{\scriptstyle
A \subset V(X) \atop A\not= \emptyset, V(X)} 
{|\partial A| \over \min (|A|, |V(X) - A|)},$$
where, for a subset
$A$ of $V(X)$, $\partial A$ is the set of edges connecting
a vertex in $A$ to one not in $A$. (See Figure 2.)

\vskip 18pt
\centerline{\psfig{figure=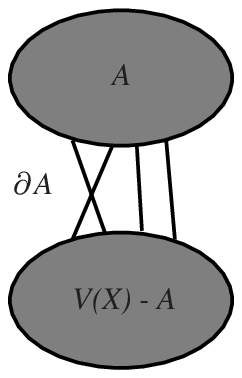}}
\vskip 18pt
\centerline{Figure 2.}

The Cheeger constant of a Riemannian manifold $M$ is defined similarly.
One considers all possible codimension one submanifolds $F$ that separate
$M$ into two manifolds $M_1$ and $M_2$. Then the {\sl Cheeger
constant} $h(M)$ is
$$h(M) = \inf_F {{\rm Area}(F) \over \min({\rm Volume}(M_1), {\rm
Volume}(M_2))}.$$

For a finite graph $X$, the smallest non-zero eigenvalue
of the Laplacian is denoted $\lambda_1(X)$. 
It is given by the formula
$$\inf \left \{ {||df||^2 \over ||f||^2} : f \hbox{ is a real-valued
function on } V(X) \hbox{ such that } \sum_{v \in V(X)} f(v) =
0 \right \},$$
where $$
\eqalign{
||f||^2 &= \sum_{v \in V(X)} |f(v)|^2 \cr
||df||^2 &= \sum_{e \in E(X)} |f(e_+) - f(e_-)|^2.\cr}$$
Here, we have picked an arbitrary orientation on each
edge $e$ of $X$, and have denoted its target and source vertex
by $e_+$ and $e_-$. It is clear that the formula for
$\lambda_1(X)$ is independent of this choice of orientation.

Similarly, if $M$ is a closed Riemannian manifold,
then $\lambda_1(M)$, the smallest non-zero eigenvalue
of the Laplacian, is given by 
$$\lambda_1(M) = \inf \left \{ {\int_M |df|^2 \over \int_M |f|^2} 
: f \in C^\infty(M), f \not\equiv 0 \hbox{ and } \int_M f = 0 \right \}.$$

The above formula for $\lambda_1(X)$
can be interpreted as follows. A subset $A$ of $V(X)$,
as in the definition of the Cheeger constant $h(X)$,
determines (up to a multiplicative factor) 
a step function $f \colon V(X) \rightarrow
{\Bbb R}$, with a single positive value on $A$
and a single negative value on $V(X) - A$, such that
the sum of $f(v)$ over all vertices $v$ in $V(X)$
is zero. The quantity $||df||^2/||f||^2$, which appears in
the definition of $\lambda_1(X)$, is clearly
a close relative of $|\partial A|/\min (|A|, |V(X) - A|)$,
which appears in the definition of $h(X)$.
However, the formula for $\lambda_1(X)$
considers all functions $f \colon V(X) \rightarrow {\Bbb R}$
that sum to zero, rather than just two-valued step functions.
In fact, the Cheeger constant and first eigenvalue
of the Laplacian are known to be related by
results of Cheeger [8], Buser [5], Brooks [4], Alon [1], Dodziuk [10],
Tanner [56], Alon and Milman [2].
A consequence of this relationship is the
following result, which can be found in [35].
This gives many equivalent definitions of
Property $(\tau)$, which was first defined by Lubotzky
and Zimmer [37]. This should be compared
with Theorem 1.1.

\noindent {\bf Theorem.} [35] {\sl Let $G$ be a finitely
generated group, generated by a finite symmetric set of
generators $S$. Let $\{ G_i \}$ be a collection of finite index normal
subgroups. Then the following conditions are equivalent,
are independent of the choice of $S$, and are known as Property
$(\tau)$ for $G$ with respect to $\{ G_i \}$:
\item{1.} There exists $\epsilon_1 > 0$ such that if
$\rho \colon G \rightarrow Aut(H)$ is a non-trivial
unitary irreducible representation of $G$ whose kernel contains
$G_i$ for some $i$, then for every $v \in H$ with $||v|| = 1$,
there exists an $s \in S$ such that $||\rho(s)v -v || \geq
\epsilon_1$.
\item{2.} There exists $\epsilon_2 >0$ such that all the
Cayley graphs $X_i$ of $G/G_i$ with respect to $S$
are $([G:G_i],|S|,\epsilon_2)$-expanders.
\item{3.} There exists $\epsilon_3 >0$ such that $h(X_i) \geq \epsilon_3$.
\item{4.} There exists $\epsilon_4 > 0$ such that $\lambda_1(X_i)
\geq \epsilon_4$.

\noindent If, in addition, $G = \pi_1(M)$ for some
compact Riemannian manifold $M$, and $M_i$ are the finite
sheeted covers corresponding to $G_i$, then the above conditions
are equivalent to each of the following:
\item{5.} There exists $\epsilon_5 > 0$ such that $h(M_i) \geq \epsilon_5$.
\item{6.} There exists $\epsilon_6 > 0$ such that 
$\lambda_1(M_i) \geq \epsilon_6$.

}

Property $(\tau)$ is usually defined, as above, for
finite index {\sl normal} subgroups of a finitely
generated group $G$. But, one can extend it
to collections $\{ G_i \}$ of finite index subgroups,
some of which may not be normal, by taking
any of (2) - (6) above as the definition, and
where $X_i$ now denotes the Schreier coset diagram
of $G/G_i$ with respect to $S$.

Thus, the failure of Property $(\tau)$ asserts that
$\lambda_1(X_i)$, $h(X_i)$, $\lambda_1(M_i)$ and $h(M_i)$ all
have subsequences that tend to zero. The force of Theorem 1.1 is that if the
convergence is fast enough (when measured by $[G:G_i]$
and $[N(G_i):G_i]$), then we can deduce positive
virtual $b_1$.

We now embark on the proof of Theorem 1.1.
Given a finite graph $X$ and two disjoint subsets $V_1$ and $V_2$
of $V(X)$, let $e(V_1,V_2)$ denote the number of edges
that have an endpoint in $V_1$ and an endpoint
in $V_2$. 

\noindent {\bf Lemma 2.1.} {\sl Let $X$ be a Cayley graph of
a finite group. Let $A$ be a non-empty subset of $V(X)$ such that
$|\partial A|/|A| = h(X)$. Then $|A| > |V(X)|/4$.}

\noindent {\sl Proof.} Consider a
smallest non-empty set of vertices $A$ satisfying
$|\partial A|/|A| = h(X)$. Suppose that
$|A| \leq |V(X)|/4$. Let $B$ be the image of $A$ 
under the left action of an arbitrary group element. Then
$$\eqalign{
&|\partial A| + |\partial B| - |\partial (A \cap B)| \cr
&\quad =e(A,A^c) + e(B,B^c) - e(A \cap B, A^c \cup B^c) \cr
&\quad =e(A \cap B, A^c \cap B^c) +e(A \cap B, B \backslash A)
+e(A \backslash B, A^c \cap B^c) + e(A \backslash B, B \backslash A)\cr
&\quad\quad + e(A \cap B, A^c \cap B^c) +e(A \cap B, A \backslash B) +
e(B \backslash A, A^c \cap B^c) + e(B \backslash A, A \backslash B) \cr
&\quad\quad- e(A \cap B, A \backslash B) - e(A \cap B, B \backslash A) 
- e(A \cap B, A^c \cap B^c) \cr}$$
\vfill\eject
$$\eqalign{
&\quad =e(A \cup B, A^c \cap B^c) + 2 e(A \backslash B, B \backslash A) \cr
&\quad =|\partial (A \cup B)| + 2 e(A \backslash B, B \backslash A).}$$
Now, $|A \cap B| \leq |A|$, and so by the minimality
of $A$,
$$|\partial (A \cap B)| \geq h(X) |A \cap B|,$$
with equality if and only if $A \cap B = \emptyset$
or $A = B$. Therefore, we deduce
that
$$\eqalign{|\partial (A \cup B)|
&= |\partial A| + |\partial B| - |\partial (A \cap B)| - 2 
e(A \backslash B, B \backslash A) \cr
&\leq h(X) (|A| + |B| - |A \cap B|) \cr
&=h(X) |A \cup B|.}$$
This must be an equality, since $A \cup B$ has size at 
most $|V(X)|/2$. Hence, we deduce that either $A$
equals $B$ or they are disjoint, and in the latter case,
there can be no edges joining $A$ to $B$. Hence, the
images of $A$ under the left action of the group form a disjoint union of
copies of $A$ with no edges between them. This implies
the graph is disconnected, which is impossible. 
$\square$

\noindent {\bf Lemma 2.2.} {\sl Let $X$ be a Cayley graph of
a finite group. Let $A$ be a non-empty subset of $V(X)$ such that
$|\partial A|/|A| = h(X)$ and $|A|\leq |V(X)|/2$.
Then the subgraphs induced by $A$ and its complement
are connected.}

\noindent {\sl Proof.}
Suppose that the subgraph induced by $A$ is not connected. Let $A_1$ be the
vertices of one component, and let $A_2$ be the
remaining vertices of $A$. Then,
$$|\partial A_1| + |\partial A_2|
= |\partial A| = h(X) |A| = h(X) |A_1| + h(X) |A_2|.$$
But, $|\partial A_i| \geq h(X) |A_i|$ for each $i$,
and hence, this must be an equality. However some $A_i$
has at most $|V(X)|/4$ vertices, which contradicts Lemma 2.1.

Now suppose that the subgraph induced by the complement of $A$ is not 
connected. Let $B_1$ be the vertices of one component, and let $B_2$
the remaining vertices of the complement. Now, $|B_1|$
and $|B_2|$ must each be less than $|V(X)|/2$.
For, if $|B_1| \geq |V(X)|/2$, say, then we may add $B_2$
to $A$ to obtain a collection of at most $|V(X)|/2$
vertices, having more vertices than $|A|$, but smaller
boundary, contradicting the definition of $h(X)$.
If we apply the definition of $h(X)$ to $B_1$ and $B_2$,
we obtain,
$$\eqalign{
e(A,B_1) &\geq h(X) |B_1| \cr
e(A,B_2) &\geq h(X) |B_2|. \cr}$$
Hence,
$$|A| = (e(A,B_1) + e(A,B_2))/h(X) \geq |B_1| + |B_2|.$$
But, $|A| + |B_1| + |B_2| = |V(X)|$,
which implies that $|A| \geq |V(X)|/2$. Thus, $|A| = |V(X)|/2$
and each of the above inequalities are equalities.
But, some $B_i$ has at most $|V(X)|/4$ vertices, which
by Lemma 2.1, implies that $e(A,B_i) > h(X) |B_i|$,
a contradiction. $\square$

\noindent {\bf Lemma 2.3.} {\sl Let $S$ and $\overline S$ be
finite sets of generators for a group $G$. Let
$\{ G_i \}$ be a collection of finite index
subgroups of $G$. Let $X_i$ and $\overline X_i$ be
the Schreier coset graphs of $G/G_i$ with respect to $S$ and
$\overline S$. Then, there is a constant $C \geq 1$,
depending only on $S$ and $\overline S$, such that
$$C^{-1} h(X_i) \leq h(\overline X_i) \leq
C \ h(X_i).$$}

\noindent {\sl Proof.} By first changing $S$ to $S \cup
\overline S$ and then $S \cup \overline S$ to $\overline S$, it
suffices to consider the case where $S \subset \overline S$.
Then $h(X_i) \leq h(\overline X_i)$, establishing
one of the inequalities. Let $S^{(n)}$ be the set
of words in $S \cup S^{-1}$ with length at most $n$. Then $\overline S
\subset S^{(n)}$ for some $n$. Since $S^{(n)} \subset
(\dots(S^{(2)})^{(2)}\dots)^{(2)}$, it suffices to
establish the lemma in the case where $\overline S = S^{(2)}$.
Let $A$ be a subset of $V(X_i)$ such that $h(X_i)
= |\partial A|/|A|$ and $|A| \leq |V(X_i)|/2$. 
Consider the corresponding subset $\overline A$ of $V(\overline X_i)$.
An edge $e$ of $\partial \overline A$ joins a vertex
in $\overline A$ to one not in $\overline A$. These two vertices
differ by an element of $\overline S$, which is
either an element of $S$, or product of two
elements of $S$. In the first case, we obtain a
corresponding edge of $\partial A$. In the second case,
we obtain two edges $e_1$ and $e_2$ of $X_i$,
one of which ($e_1$, say) must be in $\partial A$.
The other, $e_2$, must be attached to $e_1$, and so
there are at most $4|S| - 2$ choices for $e_2$, given $e_1$.
Hence, $|\partial \overline A| \leq (4|S|-1) |\partial A|$.
So, 
$$h(\overline X_i) \leq (4|S|-1)  h(X_i),$$
which proves the lemma. $\square$

Recall that any
finitely presented group has a finite triangular
presentation. This is most easily seen by viewing
the group as the fundamental group of a finite 2-complex
with a single vertex.
One subdivides all faces with more than three edges into triangles,
without adding any vertices. One way to deal with
faces with two or fewer edges is to add to the complex a 1-cell
$g$ and a 2-cell attached along the word $g^2 g^{-1}$.
Then we enlarge any face with two or fewer edges
to triangles by adding one or more letters $g$ to its
boundary. The resulting 2-complex
has the same fundamental group and specifies a triangular
presentation.

A {\sl meta-cocycle} for a graph is 
a 1-cochain for which any closed path of length at most three
evaluates to zero. If the graph
is the 1-skeleton of a 2-complex $K$ where
every 2-cell is a triangle, a meta-cocycle determines a
genuine cocycle on $K$.

\noindent {\bf Lemma 2.4.} {\sl Let $X$ be the Cayley
graph of a finite group $G$. Suppose that 
$h(X) < \sqrt{2/(3|V(X)|)}$. Then $X$ admits
a meta-cocycle that is not a co-boundary.}

\noindent {\sl Proof.} Let $A$ be a 
non-empty subset of $V(X)$ such that
$|\partial A|/|A| = h(X)$ and $|A|\leq |V(X)|/2$.

\noindent {\sl Claim.} There are
elements $g_1, \dots, g_4$ of $G$ such that
$g_i (\partial A) \cap g_j (\partial A) = \emptyset$ if $i \not= j$.

Since $X$ is a Cayley graph, each edge is oriented. 
Let $C$ denote initial vertices of $\partial A$.
As $V(X)$ is in one-one correspondence with $G$, we view
$C$ as a set of elements of $G$.
Suppose that $g (\partial A) \cap h (\partial A) \not= \emptyset$ for
some $g$ and
$h$ in $G$. Then $g C \cap h C \not= \emptyset$, and so
$h^{-1} g \in C C^{-1}$. Suppose that the claim were not true.
Then, each 4-tuple $(g_1, \dots, g_4) \in G^4$
would have $g_i C \cap g_j C \not= \emptyset$ for
some $i \not= j$. For $1 \leq i < j \leq 4$, 
define
$$\eqalign{
p_{ij} \colon G^4 &\rightarrow G \cr
(g_1, \dots, g_4) &\mapsto g_j^{-1} g_i.}$$
Then, the sets $p_{ij}^{-1}(C C^{-1})$ cover
$G^4$. Each set has size $|G|^3 |CC^{-1}|$, and so
$$ |G|^4 \leq \left ( {4 \atop 2} \right ) |G|^3 |C|^2,$$
which implies that
$$|G| \leq 6 |C|^2 \leq 6 |\partial A|^2
= 6 (h(X) |A|)^2 < 6 \left (
\sqrt{{2 \over 3 |V(X)|}} 
{|V(X)| \over 2} \right )^2 = |G|,$$
which is a contradiction. This proves the claim.

Since $|A| > |V(X)|/4$, by Lemma 2.1, and $i$ and $j$ range from 1 to 4, 
we have $g_i A \cap g_j A \not= \emptyset$
for some $i \not= j$. Since $|A| \leq |V(X)|/2$,
$g_i(V(X) - A) \cap g_j(V(X) - A) \not= \emptyset$.
Let $c'$ be the coboundary of the characteristic function of $g_iA$,
considered as a 0-cochain on $G$. 
This has support precisely $g_i(\partial A)$. Let $c$ be
the cochain that agrees with $c'$ on those edges with
both endpoints in $g_j A$, and is zero elsewhere.
(See Figure 3.)

\vskip 6pt
\centerline{\psfig{figure=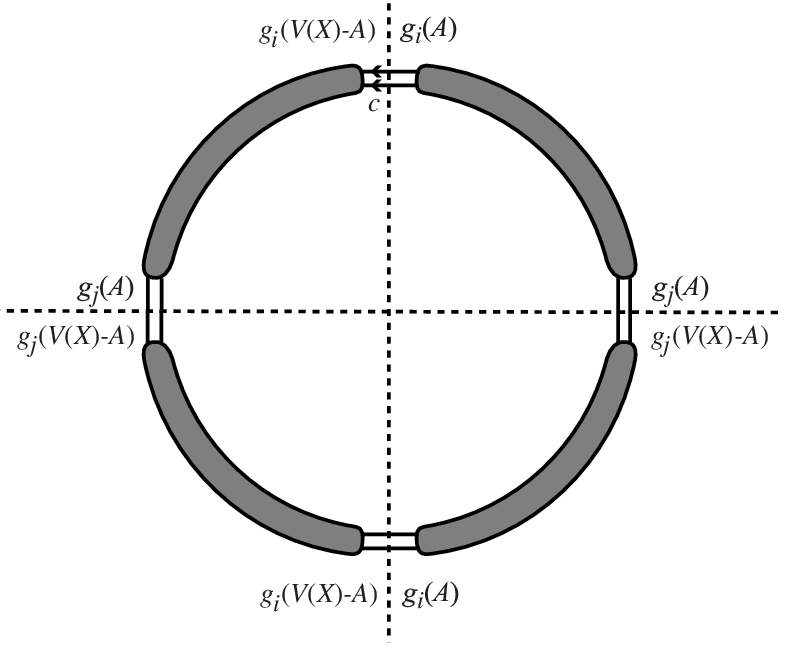}}
\vskip 18pt
\centerline{Figure 3.}

\noindent {\sl Claim.} $c$ is a meta-cocycle.

Loops of length one or two trivially evaluate to zero.
Consider a loop $\ell$ of length three in $X$. We must show that
it evaluates to zero under $c$. Since $c'$ is a coboundary,
the evaluation of $\ell$ under $c'$ is zero. As $c'$
only takes values in $\{ -1, 0, 1 \}$, $\ell$ must run
over exactly two edges in the support of $c'$. These edges
lie in $g_i(\partial A)$. Neither of these edges can lie in
$g_j(\partial A)$. Hence, either all or none of the vertices
of $\ell$ lie in $g_j(A)$. In the former case, $c$ agrees
with $c'$ on $\ell$, and hence evaluates to zero. In the
latter case, $c$ is zero on $\ell$, by definition.

\noindent {\sl Claim.} There are edges $e_1$ and $e_2$ in $g_i(\partial
A)$,
such that the endpoints of $e_1$ both lie in $g_j(A)$,
whereas neither endpoint of $e_2$ lies in $g_j(A)$.

There exists a vertex in $g_i(A) \cap g_j(A)$. Now,
$g_i(A) \not= g_j(A)$, since $g_i(\partial A) \cap
g_j(\partial A) = \emptyset$. So there exists a vertex
in $g_j(A) - g_i(A)$. By Lemma 2.2, we may pick a path in the subgraph induced
by $g_j(A)$ between these two vertices. This must run
across an edge in $g_i(\partial A)$, which we may take
to be $e_1$. To find $e_2$, apply a similar argument
joining vertices in $g_i(V(X) - A) \cap g_j(V(X) - A)$
and $g_i(A) - g_j(A)$ by a path in $g_j(V(X) - A)$.

\noindent {\sl Claim.} $c$ is not a co-boundary.

It suffices to give a loop in $X$ which has non-zero
evaluation. Now, the subgraphs induced by $g_i(A)$ and $g_i(V(X) - A)$ 
are both connected, by Lemma 2.2, and so there are arcs $\alpha_1$ and $\alpha_2$
in these subgraphs, joining
the endpoints of $e_1$ and $e_2$. The resulting loop
$e_1 \cup \alpha_1 \cup e_2 \cup \alpha_2$ has non-zero
evaluation under $c$, since its intersection with the support of
$c$ is $e_1$. $\square$

\noindent {\sl Proof of Theorem 1.2.} Let $K$ be the Cayley
2-complex arising from the finite triangular presentation,
and let $X$ be its 1-skeleton. If $h(X) < \sqrt{2/(3|V(X)|)}$,
then by Lemma 2.4, $X$ admits a meta-cocycle that is not
a co-boundary. Hence, $H^1(K)$ is infinite, but $\pi_1(K)$
is trivial, which is a contradiction. $\square$

At several points in the proof of Theorem 1.1, we will
need to consider the following situation. Let $G$ be a
finitely generated group, and let $H \triangleleft K \leq G$
be finite index subgroups. Let $X(G/H)$ and $X(G/K)$ be
the Schreier coset diagrams for $G/H$ and $G/K$ with respect
to some finite generating set $S$ for $G$. Now the group
$K/H$ acts on the right cosets $G/H$ by left multiplication,
since $H$ is normal in $K$, and hence $K/H$ 
acts freely on $X(G/H)$. The quotient of $X(G/H)$ by
this action is $X(G/K)$, and the quotient map $X(G/H) \rightarrow
X(G/K)$ is a covering map. Pick a maximal tree $T$ in
$X(G/K)$. Its inverse image in $X(G/H)$ is a forest $F$.
If we collapse each component of this forest to a
single vertex, we obtain a Cayley graph $X(K/H)$
for $K/H$. (See Figure 4.)

\vskip 18pt
\centerline{\psfig{figure=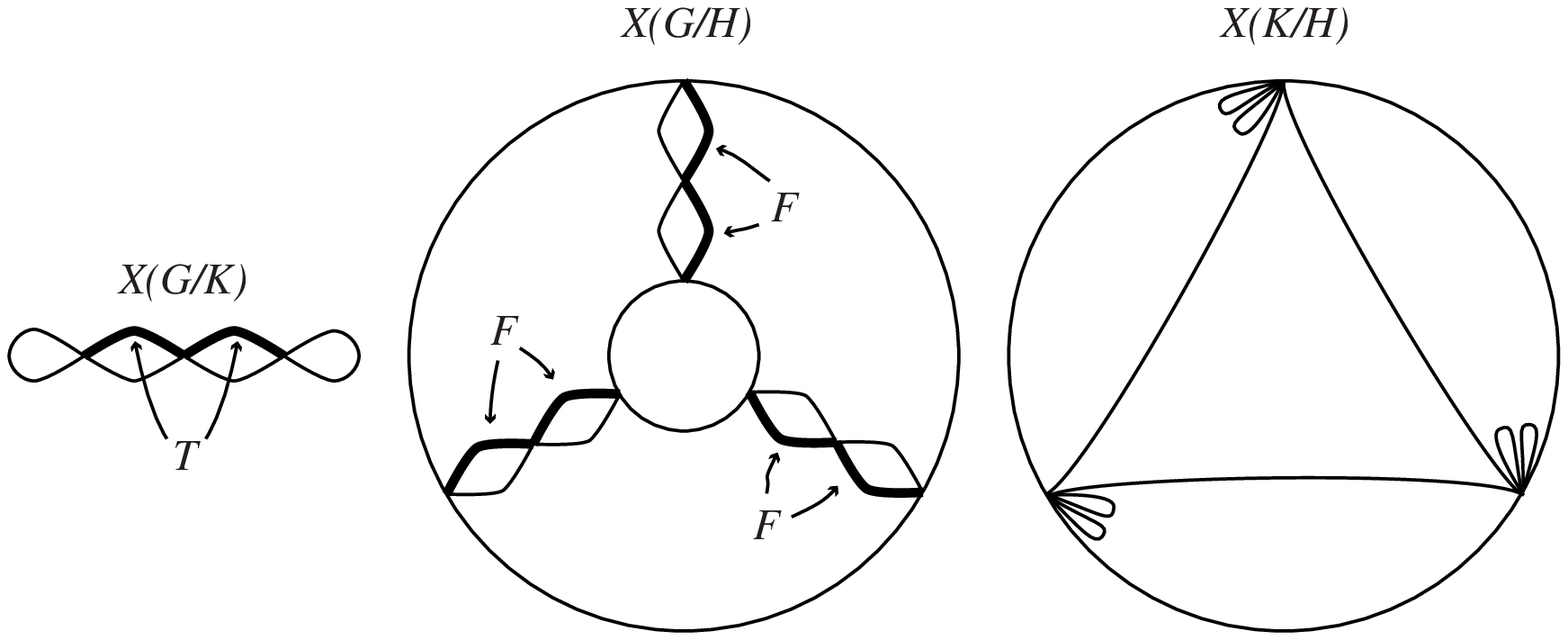,width=4in}}
\vskip 18pt
\centerline{Figure 4.}

\noindent {\bf Lemma 2.5.} {\sl With the above notation,
$$\lambda_1(X(G/H)) \leq \lambda_1(X(K/H)) \ [G:K]^{-1}.$$
Provided $h(X(G/H)) \leq [G:K]^{-1}/2$,
$$h(X(K/H)) \leq 4 \ h(X(G/H)) \ [G:K]^2 \ |S|.$$}

\noindent {\sl Proof.} Let $p \colon X(G/H) \rightarrow
X(K/H)$ be the quotient map. Any real-valued function $f$ 
on $V(X(K/H))$ induces a function
$f \circ p = \tilde f$ on  $V(X(G/H))$.
If $f$ sums to zero, so does $\tilde f$. This is
because every component of $F$ has the same number
of vertices, $[G:K]$. Also,
$$\eqalign{
||d \tilde f||^2 & = ||df||^2 \cr
|| \tilde f ||^2 & = [G:K] \ ||f||^2.}$$
The first inequality of the lemma follows immediately.

We now prove the second inequality of the lemma.
Let $A$ be a non-empty subset of $V(X(G/H))$
such that $|A| \leq |V(X(G/H))|/2$ and $h(X(G/H)) 
=|\partial A|/ |A|$. Then $|A| \geq 1/h(X(G/H))$.
We use $A$ to construct a subset $B$ of $V(X(K/H))$.
If some component of the forest $F$ has all its vertices
in $A$ (respectively, all its vertices in $A^c$),
we place the corresponding vertex of $X(K/H)$ in
$B$ (respectively, $B^c$).
We place the remaining vertices of $X(K/H)$
in $B$ or $B^c$, so that
$|B|$ is the largest integer less
than or equal to $|A|/[G:K]$. So,
$|B| \leq |V(X(K/H))|/2$. Also,
$$|B| \geq {|A| \over [G:K]} - 1
\geq {|A| \over 2 [G:K]}.$$
The second inequality holds, since 
$${|A| \over 2 [G:K]} \geq {1/h(X(G/H)) \over 2 [G:K]} \geq 1,$$
by assumption.

We now wish to estimate $|\partial B|$.
Each edge $e'$ in $\partial B$
comes from a single edge $e$ in $X(G/H)$ with
endpoints in distinct components of $F$. If $e$  lies
in $\partial A$, we have a contribution to $\partial A$.
If $e$ does not lie in $\partial A$, at least
one of the components of $F$ at its endpoints
does not lie wholly in $A$ or wholly in $A^c$.
Hence, this component of $F$ contains an edge $e''$
in $\partial A$. The number of edges $e'$
that can count towards $e''$ is at most the number of
edges emanating from a component of $F$, which
is less than $2 |S| [G:K]$. So,
$|\partial B| \leq 2 |S| [G:K] |\partial A|$. Hence,
$$h(X(K/H)) \leq {|\partial B| \over |B|}
\leq { 2 |S| [G:K] |\partial A|
\over |A| / (2[G:K])}
= 4 h(X(G/H)) \ [G:K]^2 \ |S|,$$ 
which proves the lemma. $\square$

\noindent {\sl Proof of Theorem 1.1.} $(1) \Rightarrow (2)$:
If some $G_i$ has infinite abelianisation, then it
admits a non-trivial homomorphism $G_i \rightarrow {\Bbb Z}$. Hence,
so does any finite index subgroup of $G_i$. There are
infinitely many finite index subgroups, for example arising
from the kernel of $G_i \rightarrow {\Bbb Z} \buildrel
{\rm mod} \ n \over \longrightarrow {\Bbb Z}_n$.

$(2) \Rightarrow (1)$: Trivial.

$(1) \Rightarrow (3)$: Suppose that $G_i$ has infinite
abelianisation, and let $\phi \colon G_i \rightarrow {\Bbb Z}$
be a surjective homomorphism. Let $G_i^n$ be the kernel
of the homomorphism $G_i \buildrel \phi \over \longrightarrow
{\Bbb Z} \buildrel {\rm mod} \ n \over \longrightarrow
{\Bbb Z}_n$. Then $G_i^n \triangleleft G_i \leq G$ are
finite index subgroups. Now, $G_i/G_i^n$ is the
cyclic group of order $n$. Hence, each vertex
$v$ of $X(G_i/G_i^n)$ is labelled with an integer
$\psi(v)$ mod $n$. (See Figure 1.)

Let $T$ be the maximal tree in $X(G/G_i)$ defined
before Lemma 2.5 (where $K = G_i$). The edges of
$X(G/G_i)$ not in $T$ give a set of generators
for $G_i$. Let $N$ be the maximal absolute value of
$\phi(g)$, as $g$ varies over these generators.
Then, each edge of $X(G_i/G_i^n)$ joins vertices
whose labels differ by at most $N$.

We will now find an upper bound for $\lambda_1(X(G_i/G_i^n))$.
Define
$$\eqalign{
f \colon V(X(G_i/G_i^n)) &\rightarrow {\Bbb R} \cr
v &\mapsto \sin (2 \pi \psi(v) / n).}$$
Note that 
$$\sum_{v \in V(X(G_i/G_i^n))} f(v) = 0.$$
Now each edge $e$ of $X(G_i/G_i^n)$, oriented
from $e_-$ to $e_+$, satisfies
$$\eqalign{
|f(e_+) - f(e_-)| 
& = \vert \sin ( 2 \pi \psi(e_+) / n)
- \sin ( 2 \pi \psi(e_-)  / n) \vert \cr
&= 2 \vert \sin (\pi ( \psi(e_+) - \psi(e_-)) / n)
\cos (\pi (\psi(e_+)+ \psi(e_-)) / n ) \vert \cr
& \leq 2 \pi N/n .}$$
Hence,
$$
||df||^2 = \sum_{e \in E(X(G_i/G_i^n))} |f(e_+) - f(e_-)|^2
\leq {[G:G_i] \over n} (4 \pi^2 N^2 |S|),$$
since $X(G_i/G_i^n)$ has at most $[G:G_i] n |S|$ edges.
Also,
$$
||f||^2 = \sum_{v \in V(X(G_i/G_i^n))} |f(v)|^2
= \sum_{m=1}^n \sin^2(2 \pi m /n) \geq n/8,$$
provided $n > 2$. The last inequality holds because at least
half the values of $2\pi m/n$ lie in the intervals
$[\pi/6,5\pi/6]$ and $[7\pi/6, 11\pi/6]$. So,
$$
\lambda_1(X(G_i/G_i^n)) \leq {||df||^2 \over ||f||^2}
\leq {[G:G_i] \over n^2} (32 \pi^2 N^2 |S|).$$
Now apply Lemma 2.5 to deduce that
$$
\lambda_1(X(G/G_i^n)) \leq {32 \pi^2 N^2 |S| \over n^2}.$$
Since $G_i^n$ is a normal subgroup of $G_i$,
$N(G_i^n)$ contains $G_i$. So, $[G:N(G_i^n)]$ is
bounded independently of $n$. Hence,
$\lambda_1(X(G/G_i^n)) [G:N(G_i^n)]^2 [G:G_i^n]^2$
is bounded independently of $n$, as required.

$(3) \Rightarrow (4)$: It is trivial to show that, if
$X$ is a finite connected graph, then $\lambda_1(X) \geq 1/|V(X)|^2$.
So, consider a subsequence for which
$\lambda_1(X_i) [G:N(G_i)]^2 [G:G_i]^2$ is bounded
above by some constant $k$. Then,
$$k \geq \lambda_1(X_i) [G:N(G_i)]^2 [G:G_i]^2
\geq [G:N(G_i)]^2.$$
So, $[G:N(G_i)]$ is bounded. Since $[G:G_i]$
tends to infinity, so must $[N(G_i):G_i]$. We deduce that
$$\lambda_1(X_i) [G:N(G_i)]^4 [N(G_i):G_i] 
\leq k / [N(G_i):G_i] \rightarrow 0.$$ 

$(4) \Rightarrow (6)$: This follows from the fact ([1],[10]) that if
$X$ is a finite graph and each of its vertices has valence
at most $2 |S|$, then $h(X) \leq \sqrt{4 |S| \lambda_1(X)}$.

$(1) \Rightarrow (5)$: Let $G_i^n$ be the subgroups in
the proof of $(1) \Rightarrow (3)$. Let $B$ be the
vertices of $X(G_i/G_i^n)$ with label between $1$
and $n/2$. Let $A$ be their inverse image in
$X(G/G_i^n)$. Then $|A| / |V(X(G/G_i^n))|
\rightarrow 1/2$ and $|A|\leq |V(X(G/G_i^n))|/2$, 
but $|\partial A|$ is uniformly bounded.
So, 
$$\eqalign{
h(X(G/G_i^n)) \ [G:N(G_i^n)] \ [G:G_i^n]
\leq [G:G_i] \ |V(X(G/G_i^n))| \ |\partial A| \ / \ |A|,\cr}$$
which is uniformly bounded.

$(5) \Rightarrow (6)$: This is similar to $(3) \Rightarrow (4)$.
Since $h(X) \geq 2/|V(X)|$ for a finite connected graph $X$,
the hypothesis that $h(X_i) [G:N(G_i)][G:G_i]$
is bounded implies that $[G:N(G_i)]$ is bounded.
Hence, $[N(G_i): G_i]$ tends to infinity,
and therefore $h(X_i) [G:N(G_i)]^2 [N(G_i):G_i]^{1/2}$
has zero infimum.

$(6) \Rightarrow (1)$: This is the
heart of Theorem 1.1, where we go from a geometric
hypothesis to an algebraic conclusion.

By Lemma 2.3, if (6) holds for some finite presentation of $G$,
it holds for any finite presentation. So, we use
a finite triangular presentation.

Since the quantity in (6) has zero infimum, so must
$h(X_i) [G:N(G_i)]$. Hence, we may apply
Lemma 2.5, with $H = G_i$ and $K = N(G_i)$,
to deduce that 
$$h(X(N(G_i)/G_i)) \leq 
4 h(X_i) \ [G:N(G_i)]^2 \ |S|. $$ Therefore
$h(X(N(G_i)/G_i)) [N(G_i) : G_i]^{1/2}$ tends to zero in some
subsequence. For sufficiently
large $i$, we may apply Lemma 2.4 to deduce that $X(N(G_i)/G_i)$ has
a meta-cocycle that is not a co-boundary.
This pulls back to a meta-cocycle for $X_i$.
It also is not a co-boundary. For, there is a closed
loop in $X(N(G_i)/G_i)$ that has non-zero evaluation, and we may
use this to construct a closed loop in $X_i$ with
non-zero evaluation. Since the presentation for $G$
was triangular, this gives a non-trivial cocycle
for $G_i$, which implies that it has infinite abelianisation.

Thus, we have established the equivalence of $(1)$ - $(6)$. In
particular, $(3)$ - $(6)$ do not depend on the choice
of generators $S$ for $G$.

Suppose now that $G$ is the fundamental group of some
closed orientable Riemannian manifold $M$, and let $M_i$ be
the cover of $M$ corresponding to $G_i$.
In [4], a finite set $S$ of generators for
$\pi_1(M)$ are chosen that arise from a fundamental domain.
By the argument of Lemma 2 in [4],
$h(X_i) \leq c \ h(M_i)$, for a constant $c$ independent of
$i$. Moreover, it is simple to ensure that
$h(M_i) \leq c \ h(X_i)$ and that $\lambda_1(M_i) \leq c \lambda_1(X_i)$.

$(3) \Rightarrow (7)$; $(4) \Rightarrow (8)$; 
$(5) \Leftrightarrow (9)$; $(6) \Leftrightarrow (10)$:
These follow immediately.

$(7) \Rightarrow (9)$; $(8) \Rightarrow (10)$: These
follow from Cheeger's inequality [8], which states that
$\lambda_1(M_i) \geq (h(M_i))^2/4$.
$\square$

\vskip 18pt
\centerline {\caps 3. Background on Heegaard splittings}
\vskip 6pt

We now recall various definitions relating to Heegaard
splittings. We then collate some fairly elementary information about
the behaviour of Heegaard splittings when they are lifted to
finite covers. 

Recall that a {\sl handlebody} is a compact orientable
3-manifold obtained from a 3-ball by attaching
a (possibly empty) collection of 1-handles. A {\sl compression body} is a 
connected compact orientable 3-manifold that either is
a handlebody or is obtained from $F \times I$,
where $F$ is a possibly disconnected, closed, orientable
surface, by attaching 1-handles to $F \times \{ 1 \}$.
The {\sl negative boundary} is empty in the
case of a handlebody, and $F \times \{ 0 \}$
otherwise. The {\sl positive boundary} is
the remaining boundary components. A {\sl meridian
disc} in a compression body is a properly
embedded disc, whose boundary is an essential
curve in the positive boundary.

A {\sl Heegaard splitting} of a compact orientable
3-manifold is a description of the manifold
as two compression bodies glued along their
positive boundary. The latter forms the
{\sl Heegaard surface}. Any compact orientable
3-manifold has a Heegaard splitting, and
therefore its Heegaard Euler characteristic,
defined in \S1, is well-defined. A
Heegaard splitting is {\sl reducible} if
there are meridian discs in each compression
body with equal boundary; otherwise
it is {\sl irreducible}. It is {\sl
weakly reducible} if there are meridian
discs in each compression body with
disjoint boundaries; otherwise it is
{\sl strongly irreducible}. It is a
well known result of Casson and Gordon [6]
that if a 3-manifold admits a Heegaard
splitting which is irreducible but weakly
reducible, it is Haken. Thus, for irreducible
non-Haken manifolds $M$, $\chi_-^h(M)$ and
$\chi_-^{sh}(M)$ are equal.

A Heegaard splitting arises from a handle
structure on the 3-manifold. For, any closed
orientable 3-manifold is obtained from a collection
of 0-handles, by attaching 1-handles, then
2-handles, then 3-handles. The 3-manifold
obtained after attaching all the 1-handles
is a handlebody, as is the closure of its
complement. Thus, the boundary of this submanifold
is a Heegaard surface. 

By considering more general handle structures,
Scharlemann and Thompson introduced generalised
Heegaard splittings [51].
Here, one views a compact orientable 3-manifold $M$ as built from 
a collection of 0-handles, and possibly collars on some boundary components
of $M$, by attaching a collection of 1-handles, then 
a collection of 2-handles, then 1-handles, and so on
in an alternating fashion, and ending finally in a collection of
3-handles.
If one were to halt this process just after attaching
the $j^{th}$ batch of 1-handles or 2-handles, the result would be
a 3-manifold embedded in $M$. Let $F_j$ be the boundary of
this 3-manifold, after discarding any 2-sphere components
that bound a 0-handle or 3-handle in $M$.
After a small isotopy to make them disjoint, the surfaces 
$F_{2j}$ (known as {\sl even} surfaces)
divide $M$ into 3-manifolds, for which the surfaces $F_{2j+1}$ 
(known as {\sl odd} surfaces) form Heegaard surfaces. 
We term the number of such 3-manifolds
the {\sl length} of the generalised Heegaard splitting.
Thus, a splitting of length one is a genuine Heegaard
splitting. A schematic diagram of a generalised Heegaard
splitting is shown in Figure 5.

\vskip 18pt
\centerline{\psfig{figure=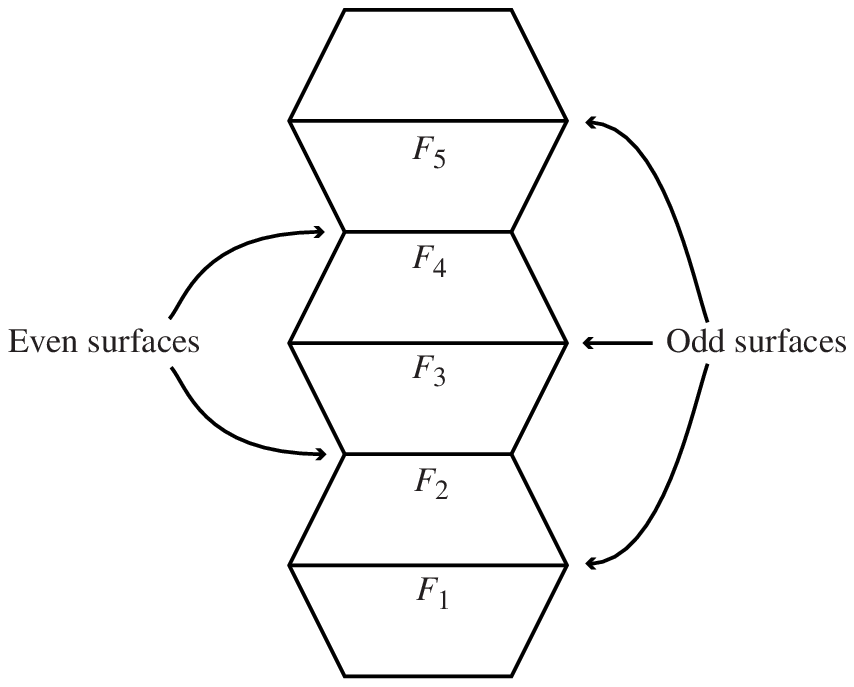,width=2.3in}}
\vskip 18pt
\centerline{Figure 5.}

Define the {\sl complexity} $c(S)$ of a closed orientable
connected surface $S$ to be $\max \{ 0, 1 -\chi(S) \}$.
Define the {\sl complexity} $c(S)$ of a closed orientable 
disconnected surface $S$ to be the sum of the complexities
of its components. The {\sl complexity} of a generalised
Heegaard splitting $\{ F_1, \dots, F_n \}$ is defined
in [51] to be 
$$\{ c(F_{2j+1}) : 1 \leq 2j+1 \leq n \},$$
where repetitions are retained. These sets are compared
lexicographically, comparing the largest integers first
and working down. This is a well-ordering. 
A decomposition of minimal complexity
is known as {\sl thin}. If $\{ F_1, \dots, F_n \}$
is a thin decomposition for $M$, we set $c_+(M)$
to be the maximal value of $c(F_{2j+1}) - 1$, as 
$2j+1$ ranges from $1$ to $n$.

It is proven in [51], using
results of Casson and Gordon [6], that for a thin generalised
Heegaard splitting of length more than one, the even surfaces 
are incompressible. In fact, from any Heegaard
surface $F$ for an irreducible 3-manifold $M$, 
one can construct a generalised Heegaard splitting
$\{ F_1, \dots, F_n \}$ with the following properties:
\item{1.} each odd surface is strongly irreducible;
\item{2.} each even surface is incompressible and has no 2-sphere
components;
\item{3.} $F_{j+1}$ and $F_j$ are not parallel for any $j$;
\item{4.} $\sum_j (-1)^j \chi(F_j) \leq -\chi(F)$;
\item{5.} $\{ c(F_{2j+1}) : 1 \leq 2j+1 \leq n \} \leq c(F)$.

\noindent This procedure is described in [51]; we summarise it
here. At each stage, starting with
the given Heegaard surface, a generalised Heegaard splitting is constructed, 
with lower complexity than the one before, and hence the procedure
is guaranteed to terminate. So, let $\{ F_1, \dots, F_n \}$
be some generalised Heegaard splitting.
If some $F_{2j+1}$ were weakly reducible,
then the modification in Rule 3 of [51] would replace $F_{2j+1}$
with three surfaces $G_1$, $G_2$ and $G_3$,
where $\chi(G_1) = \chi(G_3) = \chi(F_{2j+1})+2$
and $\chi(G_2) = \chi(F_{2j+1})+4$. Note that this
leaves $\sum_j (-1)^j \chi(F_j)$ unchanged. (See Figure 6.)
Thus, by applying
this procedure enough times, we end with
a generalised Heegaard splitting $\{ F_1, \dots, F_n \}$ 
in which each odd surface is strongly
irreducible. Hence, by Rule 5 in [51], each even surface is 
incompressible. If some even surface has a 2-sphere component,
it bounds a 3-ball since $M$ is irreducible. Pass to an
innermost such 2-sphere $S$. A component of an odd surface
forms a Heegaard surface for the 3-ball that $S$ bounds.
Remove this component and $S$. This reduces complexity
and retains the remaining conditions. If some $F_j$ and $F_{j+1}$ are
parallel, they are both discarded. This reduces
complexity and preserves the remaining conditions.
So, we end with a generalised Heegaard splitting satisfying
conditions (1) to (5).

\vskip 18pt
\centerline{\psfig{figure=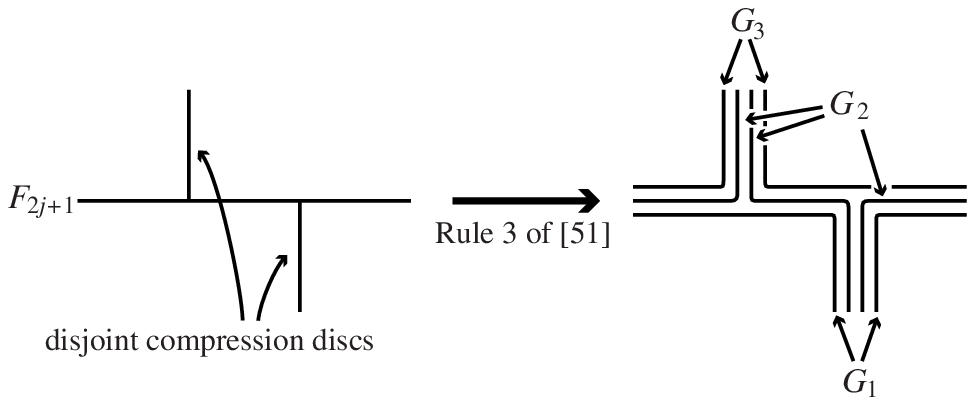}}
\vskip 18pt
\centerline{Figure 6.}

We now go on to investigate how Heegaard splittings
behave under finite covers. None of the results 
in the remainder of this section will be used later 
in this paper, but they do form a useful backdrop. 
The following, fairly well-known, observation was the
motivation for much of this paper. 

\noindent {\bf Proposition 3.1.} {\sl Let $M$ be a
compact orientable 3-manifold. Let $F$ be any Heegaard
splitting for $M$. Then, for all covers of $M$ with
sufficiently large degree, $F$ lifts to a weakly reducible
Heegaard surface.}

\noindent {\sl Proof.} Pick a maximal collection of
non-parallel meridian discs for each side of the
Heegaard surface. Suppose that there are $n_1$ discs
on one side, and $n_2$ on the other. Let $k$ be the
total number of intersection points between these discs.
Consider a $d$-fold cover $\tilde M \rightarrow M$, and
let $\tilde F$ be the inverse image of $F$. Then,
the discs lift to $dn_1$ and $dn_2$ meridian
discs either side of $\tilde F$. The total number
of intersections between these discs is $dk$. When
$(dn_1) (dn_2) > dk $, at least one pair of discs,
one on each side of $\tilde F$, are disjoint. Hence,
$\tilde F$ is weakly reducible. $\square$

However, the above fact is certainly not enough on
its own to prove the virtually Haken conjecture,
since the lifted Heegaard surface may be reducible.
In fact, the following proposition asserts that 
generically this will be the case.

\noindent {\bf Proposition 3.2.} {\sl Let $M$ be a closed orientable
Riemannian 3-manifold, and let $F$ be a Heegaard surface
for $M$. Then, for all finite covers $\tilde M$ of sufficiently
large injectivity radius, the inverse image $\tilde F$ of $F$ is
a reducible Heegaard surface.}

\vfill\eject
\noindent {\sl Proof.} If $M$ is reducible, then so is
any cover $\tilde M$, and hence so is $\tilde F$ [21]. 
Hence, we may assume that $M$ is irreducible. This implies
that $\pi_1(M)$ is not free [27].

We may take $F$ to be the boundary of a small regular neighbourhood
of a bouquet of circles $W$ embedded in $M$. This specifies
a set of generators for $\pi_1(M)$. Since $\pi_1(M)$ is
not free, there is some non-trivial word in these generators that is
trivial in $\pi_1(M)$. Pick a shortest such word. This specifies a loop $\alpha$
in $W$. Let $l$ be its length in the path metric 
on $W$. 

Pass to a cover $\tilde M \rightarrow M$ with injectivity
radius more than $l$. Let the basepoint of $M$ be the basepoint of
$W$, and assign a basepoint for $\tilde M$ 
which lies in the inverse image of the basepoint of $M$.
Consider the ball of radius $l$ in $\tilde M$ about this
basepoint. This is a 3-ball $B$. Since $\alpha$ is homotopically trivial
in $M$, it lifts to a loop $\tilde \alpha$ based at the
basepoint of $\tilde M$. As $\alpha$ represents a shortest
non-trivial word in the generators that is trivial in $\pi_1(M)$,
$\alpha$ is an embedded loop in $B$. Hence, by a theorem
of Frohman [17], $\tilde F$ is reducible. $\square$

\noindent {\bf Corollary 3.3.} {\sl The infimal Heegaard gradient
of a closed hyperbolic 3-manifold is strictly less than its
Heegaard Euler characteristic.}

\noindent {\sl Proof.} Let $F$ be a Heegaard surface for $M$
with $|\chi(F)| = \chi_-^h(M)$. The fundamental group of a hyperbolic
3-manifold is residually finite, and so we may find
finite covers with arbitrarily large injectivity radius.
By Proposition 3.2, we may find a cover $\tilde M \rightarrow M$, 
with degree $d$, say, such that the inverse image $\tilde F$ of $F$ is reducible.
However, $\tilde M$ is irreducible, and 
therefore $\tilde F$ is not minimal genus. So, $\chi_-^h(\tilde M)
< d |\chi(F)|$, which proves the corollary. $\square$

The rationale of the Heegaard gradient conjecture
is as follows. The above corollary gives that, when
passing to a finite cover of a closed hyperbolic 3-manifold, 
the Heegaard Euler characteristic 
will in general be scaled by less than the degree of the cover.
However, the conjecture proposes that it will not degenerate
too much, unless the manifold is virtually fibred.

We now give some elementary information about Heegaard
gradient. The first observation to make is that one could
very well consider any of the following quantities:
$$\eqalign{
c_1 &= \inf \{ \chi_-^h(M_i) / d_i : M_i \rightarrow M \hbox{ is
a cover with degree } d_i \} \cr
c_2 &= \inf \{ \chi_-^h(M_i) / d_i : M_i \rightarrow M \hbox{ is
a regular cover with degree } d_i \} \cr
c_3 &= \liminf \{ \chi_-^h(M_i) / d_i : M_i \rightarrow M \hbox{ is
a cover with degree } d_i \} \cr
c_4 &= \liminf \{ \chi_-^h(M_i) / d_i : M_i \rightarrow M \hbox{ is
a regular cover with degree } d_i \}.}$$
It is clear that
$$c_1 = c_2 \leq c_3 = c_4.$$
This is because any finite cover $M_i \rightarrow M$ has a further finite
cover $M_j \rightarrow M_i$ such that $M_j \rightarrow M_i
\rightarrow M$ is regular. Moreover, $\chi_-^h(M_j) / d_j
\leq \chi_-^h(M_i) / d_i$, where $d_i$ and $d_j$ are the
relevant covering degrees, since any Heegaard splitting
for $M_i$ lifts to one for $M_j$ and its Euler characteristic
is scaled by $d_j / d_i$. It is also clear that
$$c_1 = c_2 = c_3 = c_4$$
when $M$ has infinitely many covers. This is because, given any
finite cover $M_i \rightarrow M$ and an infinite sequence
of covers $M_j \rightarrow M$, we may consider the finite
covers $M_{i,j} \rightarrow M$ corresponding to the
subgroups $\pi_1(M_i) \cap \pi_1(M_j)$ of $\pi_1(M)$.
These have $\chi_-^h(M_{i,j})/d_{i,j} \leq \chi_-^h(M_i)/d_i$,
where again $d_{i,j}$ is the relevant covering degree.
When $M$ has only finitely many covers, $c_3$ and $c_4$
are both infinite, since, by convention, the $\liminf$
of a finite set of numbers is infinite.

One may also consider the quantities defined in the
same way as $c_1$, $c_2$, $c_3$ and $c_4$ but
with $\chi_-^h(M_i)$ replaced by $\chi_-^{sh}(M_i)$.
In this case, there are again obvious inequalities,
but no immediate equalities. This is because a
strongly irreducible Heegaard surface may lift
to a weakly reducible Heegaard splitting in a finite
cover. However, if $M_i$ is irreducible and non-Haken,
then $\chi_-^h(M_i) = \chi_-^{sh}(M_i)$. So,
for the purposes of proving the virtually Haken
conjecture, almost any of these quantities may
be used interchangeably.

Note that the Heegaard Euler characteristic of a closed orientable 3-manifold $M$ is always
non-negative, unless $M$ is the 3-sphere. Thus, the Heegaard
gradient of a closed orientable 3-manifold is negative if and only if its
universal cover is $S^3$. Also, its strong Heegaard gradient
is always non-negative.

There do exist manifolds with positive Heegaard gradient,
as given by the following proposition.

\noindent {\bf Proposition 3.4.} {\sl Let $M$ be a closed orientable
reducible 3-manifold other than $S^2 \times S^1$ and
${\Bbb RP}^3 \# {\Bbb RP}^3$. Then the infimal Heegaard 
gradient of $M$ is at least ${1 \over 3}$.
Moreover, when $M$ has no ${\Bbb RP}^3$ summand,
it is at least ${2 \over 3}$.}

\noindent {\sl Proof.} Pick a maximal
collection of disjoint non-parallel essential 2-spheres
in $M$. Consider a degree $d$ cover $\tilde M \rightarrow M$.
The inverse image of the 2-spheres in $M$ is a collection
${\cal S}$ of at least $d$ 2-spheres in $\tilde M$. 
If two are parallel, then between them is a copy of 
$S^2 \times I$. The only orientable manifolds that are covered by 
$S^2 \times I$ are itself and once-punctured ${\Bbb RP}^3$ [34]. 
The former is impossible in this case, since the
spheres in $M$ were not parallel and $M$ is not
$S^2 \times S^1$. We deduce that
no three spheres in ${\cal S}$ can be parallel,
since $M$ is neither $S^2 \times S^1$ nor ${\Bbb RP}^3 \# {\Bbb RP}^3$.
Also, when $M$ has no ${\Bbb RP}^3$ summand, no
two spheres of ${\cal S}$ are parallel.
If any sphere in ${\cal S}$ is not parallel to any of the
others, add in a parallel copy. Let ${\cal S}_+$ be the
resulting enlargement of ${\cal S}$. 
Let $\{ M_i \}$ be the complementary regions of ${\cal S}_+$
which are not copies of $S^2 \times I$.
Since $M$ is neither $S^2 \times S^1$ nor
${\Bbb RP}^3 \# {\Bbb RP}^3$, each component of ${\cal S}_+$
is adjacent to some $M_i$.
Let $\hat M_i$ be the closed manifold obtained from
$M_i$ by attaching 3-balls to its boundary components. Now,
by Haken's theorem [21], a Heegaard surface realizing
$\chi_-^h(\tilde M)$ may be obtained from
ones from $\hat M_i$ by attaching tubes that intersect
each component of ${\cal S}_+$ in a single
closed curve. Hence,
$$\chi_-^h(\tilde M) = \sum_i (\chi_-^h(\hat M_i) + |{\cal S_+} \cap M_i|).$$
Now, if $M_i$ has one or two boundary components, the
Heegaard Euler characteristic of $\hat M_i$ is non-negative. Hence, for each $i$,
$$\chi_-^h(\hat M_i) + 2|{\cal S}_+ \cap M_i|/3 \geq 0.$$
Summing over $i$, 
$$\chi_-^h(\tilde M) \geq \sum_i |{\cal S}_+ \cap M_i|/3
= |{\cal S}_+|/3.$$
This is at least $d/3$, and when $M$ has no ${\Bbb RP}^3$
summands, it is at least $2d/3$. $\square$

The above proposition should be compared with the Heegaard gradient conjecture.
A closed orientable reducible 3-manifold is virtually fibred 
if and only if it is $S^2 \times S^1$ or ${\Bbb RP}^3 \#
{\Bbb RP}^3$. In this case, it has zero Heegaard gradient.
But in all other cases, its Heegaard gradient is positive.

We conclude this section with an example, which illustrates
that arbitrarily large degeneration of Heegaard Euler characteristic
can occur when passing to a finite cover. More precisely,
we show that, when $\tilde M \rightarrow M$ is a degree $d$
cover, the ratio $d \chi_-^h(M)/ \chi_-^h(\tilde M)$ can
be arbitrarily large. 

Let $X$ be a once-punctured torus bundle with pseudo-Anosov
monodromy. Then $X$ admits an obvious genus three Heegaard
surface. Since $X$ is not a solid torus, its Heegaard genus is
at least two. So its Heegaard Euler characteristic is either two or four.
We now construct $M$ by Dehn filling $X$ along
a slope $s$ having distance $d$ from the slope of
the fibre. If we take $d$ to be large, $M$ will
be hyperbolic and hence its Heegaard Euler characteristic
will be two or four.
Now, the slope $s$ lifts to a simple closed curve
in the $d$-fold cyclic cover $\tilde X$ of $X$.
By Dehn filling $\tilde X$ along this slope, we obtain
a $d$-fold cover $\tilde M \rightarrow M$. Since $\tilde M$
is obtained by Dehn filling a once-punctured torus bundle,
its Heegaard Euler characteristic is at most four.
Therefore, the ratio $d \chi_-^h(M)/ \chi_-^h(\tilde M)$
can be made arbitrarily large.

\vskip 18pt
\centerline{\caps 4. The Cheeger constant and generalised Heegaard splittings}
\vskip 6pt

In this section, we establish an upper bound on the Cheeger
constant of a negatively curved Riemannian 3-manifold in terms
of its generalised Heegaard splittings. We then go on to explore
some of the consequences of this.

\noindent {\bf Theorem 4.1.} {\sl Let $M$ be a complete, finite volume
Riemannian 3-manifold. Let $\kappa < 0$ be its supremal
sectional curvature. Then
$$h(M) \leq {4 \pi c_+(M) \over |\kappa| {\rm Volume}(M)} \leq
{4 \pi \chi_-^h(M) \over |\kappa| {\rm Volume}(M)}.$$}

\noindent {\sl Proof.} Let 
$\{ F_1, \dots, F_n \}$ be a thin generalised Heegaard splitting 
for $M$. Since the even surfaces are incompressible and have
no 2-sphere components, each component that is
not boundary parallel may be isotoped to either a least area minimal
surface or the orientable double cover of an embedded non-orientable
least area minimal surface ([52], [16], [39]). Furthermore, any two such components
are either equal or disjoint after the isotopy. For odd $j$, let 
$M_j$ be the manifold between $F_{j-1}$ and $F_{j+1}$ 
(letting $M_1$ and $M_n$
be the manifolds separated off by $F_2$ and $F_{n-1}$ respectively).
There is some odd $j$ such that 
${\rm Volume}(M_1 \cup \dots \cup M_{j-2})$ and ${\rm Volume}(M_{j+2}
\cup \dots \cup M_n)$ are each at most half the volume of $M$.
(When $n=1$, we take $j$ to be 1).
Now, $F_j$ forms a Heegaard surface for $M_j$, and hence
determines a sweepout of $M_j$. In any such sweepout, there
is a surface of maximum area. The infimal possible value
for this maximum is known as a {\sl minimax} value.
According to the methods of Pitts and Rubinstein [45],
there is a minimal surface $F$ whose area is this
minimax value, and which is obtained from $F_j$ 
possibly by removing some boundary parallel
components, possibly performing some
compressions and then possibly amalgamating parallel components into
a single component. Its area, by Gauss-Bonnet, is at most
$|\kappa|^{-1} 2 \pi |\chi(F_j)| \leq 
|\kappa|^{-1} 2 \pi c_+(M)$. We may therefore
find, for any $\epsilon > 0$, a sweepout of $M_j$
whose maximal area is at most $|\kappa|^{-1} 2 \pi c_+(M)
+ \epsilon$. Some surface $F_j$ in this sweepout
divides $M$ into two pieces of equal volume. 
Hence,
$$h(M) \leq {{\rm Area}(F_j) \over (1/2) {\rm Volume}(M)}
\leq {4 \pi \, c_+(M) \over |\kappa| {\rm Volume}(M)} + {2\epsilon 
\over {\rm Volume}(M)}.$$
Since $\epsilon$ was arbitrary, the first inequality of the theorem 
is established. The second follows from the trivial observation that
$c_+(M) \leq \chi_-^h(M)$. $\square$

\noindent {\bf Theorem 1.3.} {\sl Let $M$ be a closed orientable 
3-manifold with a negatively curved Riemannian metric. Let 
$\{ M_i \rightarrow M \}$ be the finite covers of $M$.
Then the following are equivalent,
and are each equivalent to the conditions in Theorem 1.1:
\item{1.} $b_1(M_i) > 0$ for infinitely many $i$;
\item{2.} $c_+(M_i) [\pi_1 M: N(\pi_1 M_i)] $ has a bounded subsequence;
\item{3.} $c_+(M_i) [\pi_1 M: N(\pi_1 M_i)] [N(\pi_1 M_i):
\pi_1 M_i]^{-1/2}$ has zero infimum.

}

\noindent {\sl Proof.} Note that $(1)$ is simply a 
restatement of $(2)$ in Theorem 1.1.

$(1) \Rightarrow (2)$: We pick an $M_i$ with $b_1(M_i)>0$
and find a closed orientable embedded surface $S$ in $M_i$ that
represents a non-trivial element of $H_2(M_i)$.
We may assume that $S$ has no 2-sphere components.
Let $M_i^n$ be the $n$-fold cyclic cover dual to $[S]$.
We will show that $c_+(M_i^n) [\pi_1 M : N(\pi_1 M_i^n)]$
is bounded independent of $n$. Since $N(\pi_1 M_i^n)$
contains $\pi_1 M_i$, $[\pi_1 M : N(\pi_1 M_i^n)]$
is bounded above. So, it suffices to show that
$c_+(M_i^n)$ is bounded above.
Let $F$ be a Heegaard surface for $M_i - {\rm int}({\cal N}(S))$
which separates the inward-pointing boundary components
from the outward-pointing ones. We view $F$
as a surface in $M_i$. We will construct
a generalised Heegaard splitting for $M_i^n$ from the
inverse images of $F$ and $S$. Take a copy of $S$ in $M_i^n$
and let $S \times [0,1]$ be its regular neighbourhood.
Pick a handle structure on $S \times [0,{1\over 2}]$
with no 3-handles.
Its 0- and 1-handles form the first handles in the
generalised Heegaard splitting on $M_i^n$, and its
2-handles form the next handles. The remaining copies
of $F$ and $S$ form odd and even surfaces, working one way
round the cyclic cover. We end by forming a handle structure in
$S \times [{1 \over 2},1]$ using 1-, 2- and 3-handles.
By picking the handle structures on $S \times 
[0,{1\over 2}]$ and $S \times [{1 \over 2},1]$
appropriately, we can ensure that the complexity of
this splitting is less than $\max \{ c(F)+c(S)+1, 2c(S)+|S|+1 \}$,
which is independent of $n$, establishing $(2)$.

$(2) \Rightarrow (3)$: For the subsequence in $(2)$,
$[\pi_1 M: \pi_1 M_i]$ tends to infinity, but
$[\pi_1 M: N(\pi_1 M_i)]$ is bounded. So,
$[N(\pi_1 M_i): \pi_1 M_i]$ tends to infinity,
and therefore $c_+(M_i) [\pi_1 M: N(\pi_1 M_i)] [N(\pi_1 M_i):
\pi_1 M_i]^{-1/2}$ tends to zero.

$(3) \Rightarrow (1)$: By Theorem 4.1, $(3)$ implies that
$$\eqalign{
& h(M_i) [\pi_1 M: N(\pi_1 M_i)]^2 [N(\pi_1 M_i) : \pi_1 M_i]^{1/2} \cr
&\quad \leq {4 \pi c_+(M_i) [\pi_1 M: N(\pi_1 M_i)]^2 [N(\pi_1 M_i) : \pi_1 M_i]^{1/2}
\over |\kappa| {\rm Volume}(M_i)} \cr
&\quad = \left( {4 \pi \over |\kappa| {\rm Volume}(M)}
\right ) c_+(M_i) [\pi_1 M: N(\pi_1 M_i)] [N(\pi_1 M_i):
\pi_1 M_i]^{-1/2},}$$
which has zero infimum. Now apply $(10) \Rightarrow (2)$
of Theorem 1.1. $\square$

\noindent {\bf Theorem 1.5.} {\sl Let $M$ be an
orientable 3-manifold that admits a complete, negatively
curved, finite volume Riemannian metric. Let $\{ M_i \rightarrow M \}$ be a
collection of finite covers of $M$. If
$\pi_1(M)$ has Property ($\tau$) with respect to
$\{ \pi_1(M_i) \}$, then the infimal Heegaard gradient of
$\{  M_i \rightarrow M \}$ is non-zero.}

\noindent {\sl Proof.} Suppose that the Heegaard gradient of
$\{ M_i \rightarrow M \}$ is zero. Then
$$h(M_i) \leq {4 \pi c_+(M_i) \over |\kappa| {\rm Volume}(M_i)}
\leq {4 \pi \over |\kappa| {\rm Volume}(M)} {\chi_-^h(M_i)
\over [\pi_1 M : \pi_1 M_i]},$$
which has zero infimum. So, $\pi_1(M)$ fails to have
Property ($\tau$) with respect to $\{ \pi_1(M_i) \}$.
$\square$

\vfill\eject
\centerline{\caps 5. Towards the virtually Haken conjecture}
\vskip 6pt

In this section, we establish Theorem 1.7. This gives two
conditions (which conjecturally always hold for a
closed orientable hyperbolic 3-manifold) that imply the existence of
a finite Haken cover.

\noindent {\bf Lemma 5.1.} {\sl Let $F$ be a Heegaard surface
for a compact orientable 3-manifold $M$. Suppose that we can find
$d_1$ non-parallel meridian discs on one side of $F$, and $d_2$
non-parallel meridian discs on the other side of $F$, that
are all disjoint from each other. Then, either $M$
has a thin generalised Heegaard splitting with length
at least two, or $\chi_-^h(M)$ is
at most $- \chi(F) - {1 \over 3} \min \{d_1,d_2\}$.}

\noindent {\sl Proof.}
The given Heegaard splitting of $M$ has complexity
$\{ c(F) \}$. The existence of $d_1$ and
$d_2$ disjoint discs either side of $F$ allows us
to attach $d_2$ 2-handles before adding $d_1$ 1-handles.
It is trivially established, by induction on $k$, that
if a closed connected orientable surface
$F$ is compressed along a collection of $k$
pairwise non-parallel discs, then the resulting
surface $F'$ has $c(F') \leq c(F) - k/3$.
Thus, we obtain a generalised Heegaard splitting of complexity
at most $\{ c(F) - d_1/3, c(F) - d_2/3 \}$.
This is an upper bound for the complexity of any thin
generalised Heegaard splitting. If this
has length two or more, the lemma is proved.
If not, it is a Heegaard splitting, with Heegaard surface $F''$,
such that
$$\eqalign{
\chi_-^h(M) &= -\chi(F'') \leq c(F'') - 1 \leq c(F) - 1 - 
{1 \over 3} \min \{ d_1, d_2 \} \cr
&= -\chi(F) - {1 \over 3} \min \{ d_1, d_2 \}.}$$
$\square$

\noindent {\bf Theorem 1.7.} {\sl Let $M$ be a compact
orientable irreducible 3-manifold, with boundary a 
(possibly empty) collection of
tori. Let $\{ M_i \rightarrow M \}$ be a lattice
of finite regular covers of $M$. Suppose that
\item{1.} $\pi_1(M)$ fails to have Property ($\tau$) with
respect to $\{ \pi_1(M_i) \}$, and
\item{2.} $\{ M_i \rightarrow M \}$ has non-zero infimal
strong Heegaard gradient.

\noindent Then, for infinitely many $i$, $M_i$ has a
thin generalised Heegaard splitting which is not a Heegaard
splitting, and so contains a closed essential surface.
In particular, $M$ is virtually Haken.}

\vfill\eject
\noindent {\sl Proof.} Note that
if $\chi_-^{h}(M_i) < \chi_-^{sh}(M_i)$, then any
thin generalised Heegaard splitting for $M_i$
has length at least two. For if a thin generalised
Heegaard splitting has length one, it is a
minimal genus Heegaard splitting that is
strongly irreducible. Thus, we may assume
that, for all but finitely many $i$, 
$\chi_-^{h}(M_i) = \chi_-^{sh}(M_i)$. Hence,
the strong Heegaard gradient and the Heegaard
gradient of $\{ M_i \rightarrow M \}$ coincide. We work with the latter.

Suppose that $M' \rightarrow M$ is
an element of $\{ M_i \rightarrow M\}$. The set of covers in $\{ M_i \rightarrow M \}$
that factor through $M' \rightarrow M$ forms a sublattice of $\{ M_i \rightarrow M \}$.
It is clear that $M'$ has non-zero Heegaard gradient with respect 
to this sublattice, and that
$\pi_1(M')$ fails to have Property ($\tau$) with respect
to $\{ \pi_1(M') \cap \pi_1(M_i) \}$. Hence,
if we pass to a finite cover $M'$ in $\{M_i \}$, the hypotheses of
the theorem still hold. By passing to such a cover if necessary, 
we may assume that the ratio of $\chi_-^h(M)$ to the
infimal Heegaard gradient of $\{ M_i \rightarrow M \}$
is arbitrarily close to one. In particular, we may ensure that 
this ratio is at most $25/24$.

Fix a Heegaard surface $F$ of $M$ for which $- \chi(F)
=\chi_-^h(M)$. This divides $M$ into two compression bodies $H^1$ and $H^2$.
For $j=1$ and $2$, 
let ${\cal D}^j$ be a collection of disjoint compression discs 
for $F$ in $H^j$ that cut $H^j$ to a 
collar on its negative boundary or to a 3-ball.
Note that $|{\cal D}^j|$ is either $1-(\chi(F)/2)$ or
$-\chi(F)/2$, depending on whether the negative boundary 
of $H^j$ is empty or not.  (Here, we are using the hypothesis that
$\partial M$ is a possibly empty collection of tori.)
We pick a handle structure on $H^1$, by first picking
a handle structure on this collar or 3-ball, and then enlarging
it to a handle structure for $H^1$ by adding 1-handles dual to 
${\cal D}^1$. The cores of the 0-handles and 1-handles
form a graph $\Gamma$ embedded in $M$. We may assume that
the boundaries of the discs in ${\cal D}^2$ run only
over the 0-handles and 1-handles of $H^1$
and hence have a projection map to $\Gamma$. Pick a maximal
tree in $\Gamma$. If we collapse this tree to a point, a
graph $\overline \Gamma$ is obtained. The edges of 
$\overline \Gamma$ determine a set $S_-$ of generators for $\pi_1(M)$. 
Let $k$ be the maximal number of intersection points
between a disc in ${\cal D}^2$ and the union of the discs in
${\cal D}^1$,
and let $S$ be the set of words in letters of $S_-$ and their 
inverses with length at most $k$. 

The fact that $\pi_1(M)$ fails to have Property
$(\tau)$ with respect to $\{ \pi_1(M_i) \}$ implies that 
the Cayley graphs $X_i$ of $\pi_1(M)/\pi_1(M_i)$
with respect to $S$ have $\liminf_i h(X_i) = 0$.
Pass to a subsequence where $h(X_i) \rightarrow 0$.
Let $d_i$ be the order of the cover $M_i$.
Let $F_i$, $H^1_i$, $H^2_i$, $\Gamma_i$, $\overline \Gamma_i$, ${\cal D}^1_i$ and 
${\cal D}^2_i$ be the inverse images in $M_i$ of $F$, $H^1$, $H^2$, $\Gamma$, $\overline \Gamma$,
${\cal D}^1$ and ${\cal D}^2$, respectively. There is 
a graph embedding $\overline \Gamma_i \rightarrow X_i$ which
induces a bijection between their vertices.
There is also a map $\Gamma_i \rightarrow \overline \Gamma_i$ that
collapses the lifts in $\Gamma_i$ of the maximal tree.

Pick a non-empty subset $A_i$ of $V(X_i)$ such that
$|\partial A_i|/|A_i| = h(X_i)$ and $|A_i| \leq |V(X_i)|/2$.
By Lemma 2.1, $|A_i| > |V(X_i)|/4$.
Let $\overline {\cal D}^1_i$ be
the set of discs in ${\cal D}^1_i$ whose edge in $\Gamma_i$ is
mapped either to an edge in $X_i$ with an endpoint in $A_i$
or to a vertex in $A_i$.
Let $\overline {\cal D}^2_i$ be those discs
in ${\cal D}^2_i$ whose boundary maps to a loop in
$X_i$ that runs over a vertex not in $A_i$.
(See Figure 7.)

\vskip 18pt
\centerline{\psfig{figure=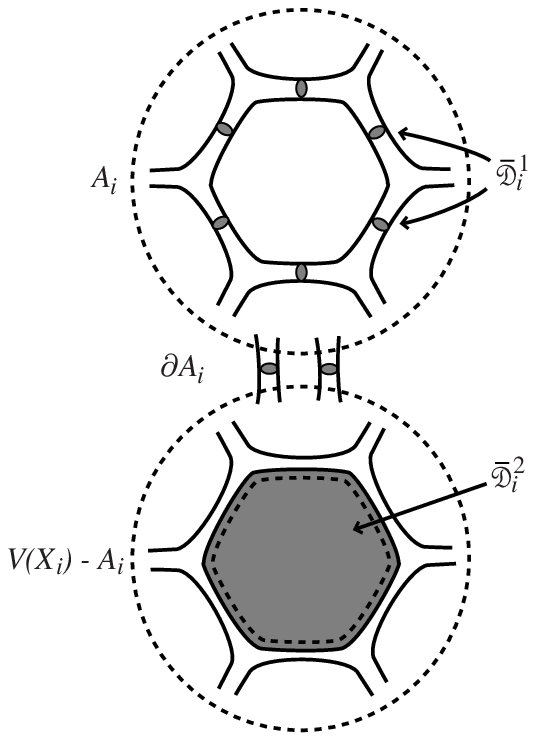}}
\vskip 18pt
\centerline{Figure 7.}

\noindent {\sl Claim.} $|\overline {\cal D}^1_i| \geq 
|{\cal D}^1| |A_i|$.

Each vertex in $A_i$ corresponds to a vertex in $\overline \Gamma_i$.
The inverse image of this vertex in $\Gamma_i$ is a tree.
There are $|{\cal D}^1|$ edges of $\Gamma_i$ that either lie
in this tree or whose initial vertices lie in the tree.
Each corresponds to a disc in $|\overline {\cal D}^1_i|$.

\noindent {\sl Claim.} $|\overline {\cal D}^2_i| \geq |{\cal D}^2| |V(X_i) - A_i|$.

This is similar to the above proof. For each disc in
${\cal D}^2$, pick a basepoint in its boundary
that maps to the vertex of $\overline \Gamma$. Then each disc in
${\cal D}^2_i$ inherits a basepoint. Each vertex
of $\overline \Gamma_i$ is in the image of $|{\cal D}^2|$
basepoints. Since $\overline {\cal D}^2_i$ contains
those discs with a basepoint in $V(X_i) - A_i$,
the inequality of the claim holds.

\noindent {\sl Claim.} The number of discs in
$\overline {\cal D}^1_i$ having non-empty intersection 
with some disc in $\overline {\cal D}^2_i$
is at most $2 |{\cal D}^1| \ |\partial A_i|$.

Suppose that two discs, $D_1$ and $D_2$, in
$\overline {\cal D}^1_i$ and $\overline {\cal D}^2_i$
intersect. The edge in $\Gamma_i$ corresponding to $D_1$ is mapped
either to an edge in $X_i$ that is attached
to some vertex $v$ in $A_i$ or to a vertex $v$ in $A_i$. 
The boundary of $D_2$ runs over
this vertex and at most $k-1$ others, one of which ($v'$, say) is
not in $A_i$. Since $S$ includes all words of length at
most $k$ in the generators $S_-$, there is an edge of
$X_i$ joining $v$ to $v'$. Thus, associated to each such pair
of discs $D_1$ and $D_2$, there is an edge in
$\partial A_i$. The number of discs $D_1$
which can correspond to this edge is at most $2|{\cal D}^1|$.
For the edge uniquely determines $v$, and there are
at most $2|{\cal D}^1|$ discs of ${\cal D}^1_i$ whose
edge in $\Gamma_i$ maps either to an edge in $X_i$ adjacent to $v$
or to $v$. This proves the claim.

Hence, if we remove all discs in $\overline {\cal D}^1_i$ that
intersect some disc in $\overline {\cal D}^2_i$, the
result is at least $|{\cal D}^1| |A_i| - 2|{\cal D}^1| |\partial A_i|$ discs.
Afterwards, none of these discs intersect. So, applying Lemma 5.1,
we deduce that either $M_i$ has a thin generalised Heegaard splitting
with length at least two, or that
$$\eqalign{
{\chi_-^h(M_i) \over d_i} &\leq
{- \chi(F_i) \over d_i} - {1 \over 3} { \min
\{ |{\cal D}^1||A_i| - 2|{\cal D}^1| |\partial A_i|, 
|{\cal D}^2| |V(X_i) - A_i| \} \over d_i }\cr
&\leq - \chi(F) + {1 \over 3} {\chi(F) \over 2} {|A_i| - 2 |\partial A_i| \over d_i}
\cr
&=\chi_-^h(M) \left (
1 - {1 \over 6} {|A_i|(1 - 2 h(X_i)) \over d_i} \right )\cr
&\leq \chi_-^h(M) \left(
{23 \over 24} + {1 \over 6} h(X_i) \right).}$$
As $i \rightarrow \infty$, the last term tends to
zero. Hence, we can find an $M_i$ such that
$\chi_-^h(M_i)/d_i < (24/25)\chi_-^h(M)$, which
is a contradiction. $\square$

\noindent {\bf Corollary 5.2.} {\sl Let $M$
be a closed orientable irreducible 3-manifold
with positive Heegaard gradient. If $\pi_1M$ is
infinite and residually finite, then it is non-amenable and
hence has exponential growth.}

\noindent {\sl Proof.} If $\pi_1(M)$ is infinite, amenable 
and residually finite, then it does not have
Property $(\tau)$ [35]. Hence, by Theorem 1.7, $M$ is
virtually Haken, and therefore satisfies the
geometrisation conjecture: it admits a decomposition
along a (possibly empty) collection of non-parallel essential tori
into Seifert fibred and hyperbolic pieces. None
of these pieces can be hyperbolic or Seifert fibred
with hyperbolic base orbifold, since $\pi_1(M)$ would
then contain a non-abelian free subgroup, contradicting
amenability. No base space can be bad or spherical, since $M$ would then
be spherical and hence have negative Heegaard
gradient. So, the base space of each Seifert
fibred piece is Euclidean. If a Seifert fibred
piece has more than one boundary component, it must be
a copy of $T^2 \times I$. As no tori in the
collection are parallel, $M$ must therefore
fibre over the circle with fibre a torus.
But it then has zero Heegaard gradient,
contrary to hypothesis. Thus, we may assume that
every Seifert fibred piece has at most one
boundary component. Hence, there are either one
or two such pieces. If there is only one, $M$
is virtually fibred, and hence has zero
Heegaard gradient, which is again a contradiction.
Thus, $M$ contains two Seifert fibred pieces
$M_1$ and $M_2$. If $\pi_1(M_1 \cap M_2)$ has
index more than two in either of $\pi_1(M_1)$
or $\pi_1(M_2)$, then $\pi_1(M)$ contains
a free non-abelian subgroup, which is
a contradiction. Hence, $\pi_1(M_1 \cap M_2)$ has
index two in each of $\pi_1(M_1)$
or $\pi_1(M_2)$, and therefore $M_1$
and $M_2$ are twisted $I$-bundles over the
Klein bottle. But $M$ is then virtually 
fibred, which again is a contradiction. $\square$

\vskip 18pt
\centerline{\caps 6. The size of minimal surfaces in negatively curved manifolds}
\vskip 6pt

In this section, we will analyse minimal surfaces in
closed negatively curved 3-manifolds. We will prove
that the diameter of such a surface is
bounded by a function of its Euler characteristic, and of the injectivity
radius and curvature of the manifold. This will have relevance
to Heegaard splittings, since it is a theorem of Pitts
and Rubinstein [45] that a strongly irreducible Heegaard surface
in a closed orientable Riemannian 3-manifold with a bumpy metric
may be isotoped either to a minimal surface or to
a double cover of a minimal non-orientable surface
with a small tube attached. In the latter case, the
tube is vertical in the $I$-bundle structure on the
regular neighbourhood of the surface.

\noindent {\bf Proposition 6.1.} {\sl 
There is some real-valued function $f$ in two variables,
with the following property. Let $M$ be a Riemannian
3-manifold, whose injectivity radius is
at least $\epsilon/2 > 0$, and whose sectional curvature
is at most $\kappa < 0$. Then the diameter of a connected
closed minimal surface $F$ in $M$ is
bounded above by $|\chi(F)| f(\kappa,\epsilon)$.}

\noindent {\sl Proof.} We give $F$ its induced Riemannian metric
as a submanifold of $M$. Since $F$ is a minimal surface, its mean
curvature is everywhere zero. This is the sum of
its principal curvatures. Hence, their product,
the extrinsic curvature of $F$, is non-positive
at all points. Now, the intrinsic curvature of $F$
at each point is the sum of its extrinsic
curvature and the ambient curvature of $M$,
and is therefore at most $\kappa <0$.

Let $F_{(0, \epsilon)}$ (respectively,
$F_{[\epsilon,\infty)}$) be the set of points in $F$ with injectivity
radius less than $\epsilon /2$ (respectively, at least $\epsilon /2$).
Here, we are referring to injectivity radius as measured
by the Riemannian metric on $F$, rather than the metric on $M$.
Given $\delta > 0$ and a metric space $X$, let $C(\delta, X)$
denote the minimal number of $\delta$-balls in $X$ required
to cover $X$. We will show that
$$\eqalignno{
C(\epsilon, F_{[\epsilon, \infty)}) &\leq 
{1 \over \cosh (|\kappa|^{1/2} \epsilon/2) - 1} |\chi(F)|, &(1)\cr
C(\epsilon + \epsilon/(2\pi) + |\kappa|^{-1/2}/2, F_{(0,\epsilon)}) 
& \leq 4 |\chi(F)|, &(2)\cr}$$
which will establish the proposition.

\vskip 18pt
\centerline{\psfig{figure=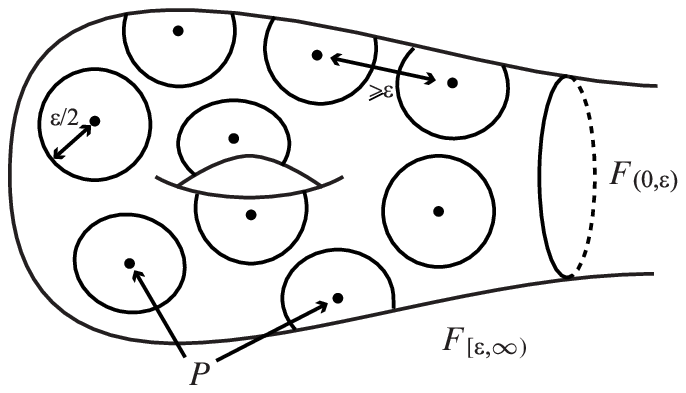}}
\vskip 18pt
\centerline{Figure 8.}

The first inequality is proved using a well-known argument
of Thurston [57]. Pick a maximal collection of points $P$
in $F_{[\epsilon, \infty)}$, no two of which lie less than
$\epsilon$ apart. (See Figure 8.) Then, the open balls of radius $\epsilon/2$ around
these points are disjoint, and each is homeomorphic to a disc,
since $P$ lies in $F_{[\epsilon, \infty)}$. By Gunther's
comparison theorem, each such disc has area at least 
$|\kappa|^{-1} 2 \pi (\cosh (|\kappa|^{1/2}\epsilon /2) - 1)$ since 
$2 \pi (\cosh (\epsilon/2) - 1)$ is the formula for the area 
of a ball of radius $\epsilon / 2$ in the hyperbolic plane.
So,
$${\rm Area}(F) \geq {\rm Area}(F_{[\epsilon, \infty)}) \geq
|P| |\kappa|^{-1} 2 \pi (\cosh (|\kappa|^{1/2}\epsilon /2) - 1).$$
By Gauss-Bonnet,
$$|\kappa| {\rm Area}(F) \leq 2 \pi |\chi(F)|.$$
By the fact that $P$ is maximal, the balls of
radius $\epsilon$ about $P$ cover $F_{[\epsilon, \infty)}$, and so
$$C(\epsilon, F_{[\epsilon, \infty)}) \leq |P| \leq
{{\rm Area}(F) \over 
|\kappa|^{-1} 2 \pi (\cosh (|\kappa|^{1/2}\epsilon /2) - 1)}
\leq 
{1 \over \cosh (|\kappa|^{1/2} \epsilon/2) - 1} |\chi(F)|,$$
which is the first inequality (1).

We now establish the second inequality (2). Pick a maximal
collection of disjoint (not necessarily simple) closed geodesics $\Gamma$ in $F$,
each with length less than $\epsilon$. Consider a geodesic
$\gamma$ in $\Gamma$. Let $\tilde F$ be the cover of $F$, corresponding to
the subgroup of $\pi_1(F)$ generated by $\gamma$. Then $\tilde F$
is homeomorphic to an open annulus, and its core geodesic $\tilde
\gamma$ is a lift of $\gamma$. Let $\tilde N$ be the set
of points in $\tilde F$ at most $\epsilon / (2\pi) + |\kappa|^{-1/2}/2$ 
from $\tilde \gamma$. (See Figure 9.)

\vskip 18pt
\centerline{\psfig{figure=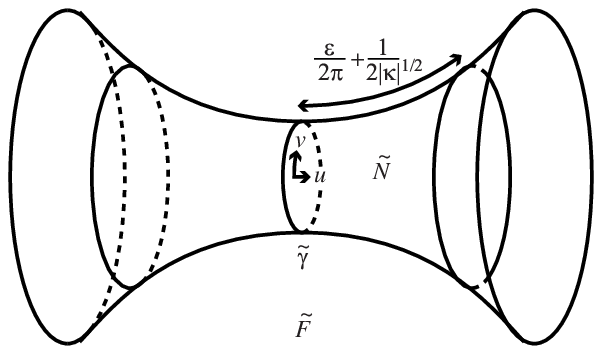}}
\vskip 18pt
\centerline{Figure 9.}

\noindent {\sl Claim 1.} Any closed curve in $\tilde F$ that is
disjoint from $\tilde N$ and that is freely homotopic to
$\tilde \gamma$ has length more than $\epsilon$.

We use the terminology of [26].
Choose orthogonal co-ordinates $(u,v)$ on $\tilde F$ so that $\tilde \gamma$
is the set $\{ u = 0 \}$ and so that the curves where $v$ is constant
are geodesics perpendicular to $\tilde \gamma$. So the metric is given by
$$ds^2 = du^2 + J^2(u,v) dv^2,$$
where $J(u,v) > 0$ and $J(0,v) = 1$. The aim is
to find a differential inequality satisfied
by the lengths of the curves $\{ u = \rm{constant} \}$. 
Thereby, we will show that, when
$u > \epsilon/(2 \pi) + |\kappa|^{-1/2}/2$,
each such curve has length more than $\epsilon$.

Let $\pm \lambda$ (where $\lambda
\geq 0$) be the principal curvatures of $\tilde F$ at a point
$(u,v)$. The intrinsic curvature $K$ at this point is the sum of the
sectional curvature of $M$ and the product of the principal
curvatures, and so
$$K \leq \kappa - \lambda^2.\eqno(3)$$
Also, by Formula 10.5.3.3 of [3], the intrinsic curvature is given by
$$K = - {1 \over J}{\partial ^2 J \over \partial u^2}.\eqno(4)$$
Let $\tilde \alpha$ be the curve $\{ u = c \}$,
where $u$ is a constant $c \geq 0$, and let $\alpha$ be its
image in $F$.
The geodesic curvature $k_g^F$ in $F$ of $\alpha$ is given by 
$$k_g^F = {1 \over J}{\partial J \over \partial u}.\eqno(5)$$
Now, $\tilde \alpha$ is freely homotopic in $\tilde F$ to
$\tilde \gamma$ which has length less than $\epsilon$ and hence maps
to a homotopically trivial curve in $M$. It therefore bounds
a mapped-in disc $D$. We may realize this disc as a minimal
surface ([11],[12],[46],[47]). 
The intrinsic curvature $K_D$ of $D$ is then at most $\kappa < 0$.
So, applying Gauss-Bonnet to $D$, we obtain
$$\int_D K_D \ dA + \int_\alpha k_g^D \ ds = 2 \pi,$$
where $k_g^D$ is the geodesic curvature of $\partial D$ in $D$.
Here, positive $k_g^D$ at a point in $\partial D$ means that
$D$ is locally convex near that point.
So, $$\int_\alpha k_g^D J \ dv = \int_\alpha k_g^D \ ds \geq 2 \pi.$$
The geodesic curvature $k_g^M$ of $\alpha$ in $M$ (which, by
definition, is non-negative) is at least $k_g^D$
and so
$$\int_\alpha k_g^M J \ dv \geq \int_\alpha k_g^D J \ dv
\geq 2\pi.\eqno(6)$$
Now, if $II$ is the second fundamental form of $F$, then
$$\tilde D_{\alpha'} \alpha' = II(\alpha',\alpha') + D_{\alpha'}
\alpha',$$
where $\tilde D$ and $D$ are the covariant derivatives in $M$ and $F$,
respectively. So, applying the triangle inequality,
$$|II(\alpha',\alpha')| \geq |\tilde D_{\alpha'} \alpha'| - 
|D_{\alpha'} \alpha'| = k_g^M J^2 - k_g^F J^2.$$
However, 
$$|II(\alpha',\alpha')| \leq \lambda |\alpha'|^2 = \lambda J^2.$$
So,
$$\lambda \geq k_g^M - k_g^F.\eqno(7)$$
So, combining (3), (4), (7) and (5), and using the universal inequality $|\kappa| +
\lambda^2 \geq 2 |\kappa|^{1/2} \lambda$:
$${1 \over J} {\partial ^2 J \over \partial u^2}
\geq |\kappa| + \lambda^2 \geq 2 |\kappa|^{1/2} \lambda 
\geq 2 |\kappa|^{1/2} k_g^M - 2 |\kappa|^{1/2} k_g^F = 
2 |\kappa|^{1/2} k_g^M - |\kappa|^{1/2} {2 \over J}{\partial J \over \partial u}.$$
Hence,
$${\partial^2 J \over \partial u^2} \geq 
2 |\kappa|^{1/2} k_g^M J - 2|\kappa|^{1/2} {\partial J \over \partial u}. \eqno(8)$$
Let $L(u)$ be the length of the curve $\alpha$.
Then
$$L(u) = \int_\alpha J \, dv.$$
So, integrating (8) with respect to $v$, and applying (6):
$${d^2 L \over d u^2} \geq 
4 |\kappa|^{1/2} \pi - 2|\kappa|^{1/2} {dL \over du}.$$
So, multiplying by $e^{2|\kappa|^{1/2} u}$ and writing $L'$ for $dL/du$:
$${d(L') \over du} e^{2|\kappa|^{1/2} u} + 2 |\kappa|^{1/2} L' 
e^{2 |\kappa|^{1/2} u} \geq 4\pi |\kappa|^{1/2} e^{2|\kappa|^{1/2} u}.$$
Integrating with respect to $u$, we obtain
$$L' e^{2 |\kappa|^{1/2} u} \geq 2\pi e^{2 |\kappa|^{1/2} u} - 2\pi,$$
and hence
$$L' \geq 2 \pi(1- e^{-2 |\kappa|^{1/2} u}).$$
Therefore,
$$L(u) \geq 2 \pi u + {\pi \over |\kappa|^{1/2}} e^{-2 |\kappa|^{1/2} u}
- {\pi \over |\kappa|^{1/2}}.$$
Now, the length of any closed curve in $\tilde F$ that is
disjoint from $\tilde N$ and that is freely homotopic to
$\tilde \gamma$ has length at least $L(\epsilon/(2\pi) + |\kappa|^{-1/2}/2)
\geq \epsilon$, which proves the claim.

Let $N_+$ be the set of points in $F$ with distance at most
$\epsilon/2 + \epsilon/(2\pi) + |\kappa|^{-1/2}/2$ from $\Gamma$. 

\vfill\eject
\noindent {\sl Claim 2.} $N_+$ contains $F_{(0,\epsilon)}$. 

Suppose that there is a point $x$ in $F_{(0,\epsilon)}$ not in $N_+$. 
There is a closed curve $\beta$ based at $x$ with length less than
$\epsilon$. Now, $\beta$ must be disjoint from $\Gamma$, since
$\Gamma$ is more that $\epsilon/2$ from $x$. Hence, $\beta$
can be freely homotoped to a closed (unbased)
geodesic $\overline \beta$ with length less than $\epsilon$
in the surface $\overline F$
obtained by cutting $F$ along $\Gamma$. Either $\overline \beta$
is disjoint from $\Gamma$, or it is a component $\gamma$ of $\Gamma$
(possibly winding several times round this component).
The former case is impossible, since $\Gamma$ is maximal.
In the latter case, $\beta$ lifts to a closed curve $\tilde \beta$ in
$\tilde F$, the cover corresponding to $\gamma$. The basepoint
of $\tilde \beta$ is at least $\epsilon/2 + \epsilon / (2\pi) +
|\kappa|^{-1/2}/2$ from $\tilde \gamma$, and hence $\tilde \beta$
misses $\tilde N$. By Claim 1, $\tilde \beta$, and hence
$\beta$, therefore has length at least $\epsilon$,
which is a contradiction, that proves the claim.

\noindent {\sl Claim 3.} The maximum number of disjoint
closed geodesics in $F$ is at most $4 |\chi(F)|$.

Let $\Gamma^+$ be a maximal collection of disjoint closed
geodesics. Let $\Gamma^+_s$ be the simple geodesics
in $\Gamma^+$ and let $\Gamma^+_{ns}$ be the non-simple
ones. Each complementary region of $\Gamma^+$ is 
an annulus, disc or M\"obius band, since any other subsurface of $F$
admits a closed geodesic that is not homotopic to a 
multiple of a boundary component. Attach the discs
to a regular neighbourhood of $\Gamma^+_{ns}$,
creating a subsurface $X$ of $F$. Each component
of $X$ has negative Euler characteristic, 
since otherwise the geodesic it contains would be simple.
Hence $|\chi(X)| \geq |X| = |\Gamma^+_{ns}|$.
Also, $|\partial X| \leq 3 |\chi(X)|$.
Each complementary region of $X$ is an annulus
or M\"obius band, that contains a single simple closed geodesic in
$\Gamma^+_s$. So, $|\Gamma^+| = |\Gamma^+_s| +
|\Gamma^+_{ns}| \leq |\partial X| + |X|
\leq 4 |\chi(X)| = 4 |\chi(F)|$,
which establishes the claim.

Hence, the number of components of $\Gamma$ is
at most $4|\chi(F)|$. Pick one point on each
component of $\Gamma$. Each point of $N_+$,
and hence $F_{(0,\epsilon)}$,
lies within $\epsilon + \epsilon/(2\pi) + |\kappa|^{-1/2}/2$
of one of these points. This proves inequality (2). $\square$

The following is an immediate application of Proposition 6.1.
It provides some evidence for the strong Heegaard gradient
conjecture.

\noindent {\bf Theorem 1.9.} {\sl Let $M$ be a closed orientable
hyperbolic 3-manifold, and let $\{ M_i \rightarrow M \}$
be the cyclic covers dual to some non-trivial element
of $H_2(M)$. Then the strong Heegaard gradient of
$\{ M_i \rightarrow M \}$ is non-zero.}

\noindent {\sl Proof.} We may perturb the hyperbolic metrics on the covering spaces
$M_i$ to bumpy metrics in the sense of [59]. Let $\kappa <0$ be
the supremum of their sectional curvatures.
Let $S$ be an incompressible surface
embedded in $M$, representing the non-trivial element of $H_2(M)$.
Let $M_S$ be the manifold obtained by cutting $M$ along $S$. Let $\delta$
be the minimal distance between the boundary
components of $M_S$. 
Let $F_i$ be a strongly irreducible Heegaard surface for $M_i$
(if one exists) with $-\chi(F_i) = \chi_-^{sh}(M_i)$.
This may isotoped 
to a minimal surface or to the boundary of a regular neighbourhood
of a non-orientable embedded minimal surface, with a small tube attached [45].
The injectivity radius of $M_i$ is at least that of $M$,
which we denote by $\epsilon/2$. By Proposition 6.1,
the diameter of $F_i$ is at most $\chi_-^{sh}(M_i) f(\kappa, \epsilon)$.
However, there are lifts
of $S$ which are at least $\lfloor i/2 \rfloor \delta$ apart.
None of these can be disjoint from $F_i$, since the
complement of $F_i$ is two handlebodies. So,
$$\chi_-^{sh}(M_i) f(\kappa,\epsilon) \geq \lfloor i/2 \rfloor \delta
\geq i\delta/3,$$
when $i \geq 2$. Rearranging this gives that
the strong Heegaard gradient is bounded below
by $\delta/(3f(\kappa,\epsilon))$. $\square$

This theorem does not generalise to cusped
hyperbolic 3-manifolds. For example, once-punctured torus
bundles have a Heegaard splitting with genus at most three.
A minimal genus splitting cannot be weakly reducible. 
For this would imply that the manifold contained
a closed essential surface [41], which does not occur [15].
Hence, if $\{ M_i \rightarrow M \}$ are the cyclic
covers dual to the fibre, then $\chi_-^{sh}(M_i)$ is at
most four. So, the strong Heegaard gradient
of $\{ M_i \rightarrow M \}$ is zero.

\noindent {\bf Corollary 1.10.} {\sl Let $M$ be a closed orientable
3-manifold that fibres over the circle with pseudo-Anosov
monodromy. Let $\{ M_i \rightarrow M \}$ be the cyclic covers dual to
the fibre. Then, for all but finitely many $i$, 
$M_i$ has an irreducible, weakly reducible, minimal genus 
Heegaard splitting.}

\noindent {\sl Proof.} By Thurston's geometrisation theorem [40],
$M$ has a hyperbolic metric. Let $S$ be the fibre. Then
$\chi_-^h(M_i) \leq 2|\chi(S)|+4$. However,
$\chi_-^{sh}(M_i)/i$ is bounded away from zero by Theorem 1.9.
So, for sufficiently large $i$, $\chi_-^{sh}(M_i) > \chi_-^h(M_i)$.
Hence, any minimal genus Heegaard splitting for $M_i$ is
weakly reducible, but necessarily irreducible. $\square$

\vfill\eject
\centerline{\caps 7. The Heegaard gradient of cyclic covers}
\vskip 6pt

The following result provides some evidence for the Heegaard gradient
conjecture.

\noindent {\bf Theorem 1.11.} {\sl 
Let $M$ be a compact orientable finite volume hyperbolic 3-manifold, and
let $\{ M_i \rightarrow M \}$ be the cyclic covers dual to some non-trivial element
$z$ of $H_2(M,\partial M)$.
Then, the infimal Heegaard gradient of $\{ M_i \rightarrow M \}$
is zero if and only if $z$ is represented by a fibre.}

\noindent {\sl Proof.} In one direction, this is trivial.
So, suppose that $z$ is not represented by a fibre.
Let $S$ be an incompressible representative for $z$
that intersects each component of $\partial M$ in a coherently
oriented (possibly empty) collection of curves.
Let $S_i$ be its inverse image in $M_i$. 

We first show that $\{\chi_-^h(M_i) \}$ has no bounded subsequence.
Suppose it does, and pass to this subsequence. Let $F_i$ be a 
Heegaard surface realizing $\chi_-^h(M_i)$, which must be
irreducible. We construct a generalised Heegaard splitting
$F^i_1, \dots, F^i_n$ from $F_i$ satisfying the five conditions in \S3.
Let $F^i$ be $F^i_1 \cup \dots \cup F^i_n$.
Because the difference in Euler characteristic between
successive odd and even surfaces is at least two,
the length of the decomposition is at most
${1 \over 2} \sum_j (-1)^j \chi(F^i_j)$, which is at most
${1 \over 2}|\chi(F_i)|$, by Condition 4 in \S3.
By assumption, $|\chi(F_i)|$ is bounded.

The argument now divides, according to whether $M$ is closed or
bounded. We start with the closed case. 
Since the even surfaces $F_{2j}^i$ are incompressible and contain
no 2-sphere components, each component may be isotoped to either a least area minimal
surface or the orientable double cover of an embedded non-orientable least area minimal
surface ([52], [16], [39]). Furthermore, after this isotopy,
any two components of the even surfaces are either equal
or disjoint.
Each component of the odd surfaces may be also be isotoped to 
a minimal surface, or to the boundary of a regular neighbourhood of an embedded minimal
surface, with a small tube attached. This follows from the argument of Pitts and Rubinstein [45], 
since the least area minimal surfaces form a barrier for the sweep-outs associated with the
odd surfaces. We now follow the argument of Theorem 1.9.
The number of components of $F^i$ is bounded, and the
Euler characteristic of each component is bounded. Hence,
by Proposition 6.1, the diameter of each component is bounded. Therefore,
the number of components of $S_i$ that $F^i$ can intersect
is bounded. However, the total number of components of $S_i$
is unbounded. Therefore, if $i$ is sufficiently large,
we may find an arbitrarily large number of copies of $S$
who components lie in the same complementary regions of $F^i$.
However, each such region is a compression body, and so
these incompressible surfaces must be parallel to 
components of the even surfaces. We deduce that two copies of
$S$ are parallel in $M_i$. Hence, $S$ is a fibre, contrary to
assumption.

We now consider the case where $M$ has non-empty boundary.
The argument in the closed case does not work, since the
injectivity radius of $M$ is zero, and so Proposition 6.1 does not apply. 
Instead, we use the theory of normal and almost normal surfaces.
Since $M$ is cusped, it has a canonical polyhedral decomposition $P$ [14],
which is an angled polyhedral decomposition in the sense of [31].
This lifts to an angled polyhedral decomposition $P_i$ of each $M_i$.
The interior angles of $P_i$ are those of $P$, and hence they have
a uniform upper bound, which is less than $\pi$, and a uniform
lower bound which is greater than zero.

We may isotope each even surface into normal form in $P_i$.
If $P_i$ were a triangulation (which it is not), then a theorem of 
Rubinstein [49] and Stocking [55] states that the odd surfaces could be placed in
almost normal form with respect to $P_i$. This was generalised in [32]
to ideal polyhedral decompositions, with a suitable generalisation
of the notion of almost normal. Two conditions on the ideal
polyhedral decomposition are required: that there is no 2-sphere
that is normal to one side, and that each face is a triangle
or bigon. The former holds, since the ideal polyhedral decomposition
is angled, and the latter holds after the faces have been
suitably subdivided. Once each $F^i_j$ is in
normal or almost normal form, it then inherits a combinatorial area, 
as in [31].
This is defined additively over each component of intersection
between $F^i$ and the polyhedra. The combinatorial area of each
such component $D$ is the sum of its exterior angles, minus 
$2\pi$ times its Euler characteristic. The area of $D$ is 
greater than or equal to zero, with equality if and only if
$D$ is the link of an ideal vertex. If the combinatorial 
area of $D$ is positive, it has a positive lower bound, $b$ say,
which applies to each of the $M_i$. By the argument of Proposition 4.3
in [31], the combinatorial area of a normal or almost normal surface
$F$ is $-2 \pi \chi(F)$. So the number of normal
or almost normal pieces of $F^i$ with positive area is
at most $2 \pi |\chi(F^i)| /b$, which is bounded.

\vfill\eject
Now, $S_i$ contains $i$ copies of $S$, which we denote by $S^i_1, \dots
S_i^i$. These are labelled so that $S^i_j$ is adjacent to $S^i_{j-1}$
and $S^i_{j+1}$, where the indexing is mod $i$. We need to consider
a closed surface, and so suppose that $S$ has non-empty boundary.
Let $k$ be $-{3 \over 2}\chi(S)+|\partial S|+1$.
Let $A^i_j$ be the annuli in $\partial M_i$ between $\partial S^i_j$ and
$\partial S^i_{j+k}$ that intersect $\partial S^i_{j+1}, \dots,
\partial S^i_{j+k-1}$. Then $S^i_j \cup A^i_j \cup S^i_{j+k}$
forms a closed surface $\hat S^i_j$ which we may push a little
into the interior of $M$. It is a theorem of Cooper and Long [9]
that $\hat S^i_j$ is essential, since $S$ is not a fibre.
Also, $\hat S^i_j$ and $\hat S^i_{j'}$
are disjoint providing $j$ and $j'$ differ by more than $k$ mod $i$.
When these surfaces are disjoint, it is clear that they cannot
be parallel. When $S_i$ is closed, let $\hat S^i_j$ be $S^i_j$.

Place $S$ in normal form in $P$. Then $S_i$ is in normal form in 
$P_i$. Suppose that each polyhedron of $P_i$ intersects $\hat S^i_j$
in at most $c$ components. Then $c$ may be taken to be independent
of $i$ and $j$. We may find $\lfloor i/(k+1) \rfloor$ disjoint
$\hat S^i_j$ in $M_i$. Now, the number of polyhedra in $P_i$ containing
positive area pieces of $F^i$ is bounded. Hence, if $i$
is sufficiently large, we may find an arbitrarily large number of
disjoint $\hat S^i_j$ with the property that the polyhedra they lie in
contain no positive area pieces of $F^i$. Since the $\hat S^i_j$
are closed, they may be isotoped off $F^i$. They then lie in the
compression bodies of the generalised Heegaard splittings, and
hence are parallel to the even surfaces. Since the total number of
components of the even surfaces is uniformly bounded,
we deduce that two disjoint $\hat S^i_j$ are parallel, which
is a contradiction. 

Thus, in both the case where $M$ is closed and the case
where $M$ has boundary,
$\chi^h_-(M_i)$ cannot have a bounded subsequence.
Hence there is a positive integer $n$, so that when $i \geq n$,
$\chi^h_-(M_i) \geq 36 |\chi(S)| + 4|S|$. This implies that the
manifold $M(i)$ obtained by cutting $M_i$ along a copy of $S$
has $\chi^h_-(M(i)) \geq 35 |\chi(S)|$. For it is not hard construct
a Heegaard surface for $M_i$ from one for $M(i)$, changing
the Euler characteristic by at most $|\chi(S)| + 4|S|$.

We wish to apply a result of Schultens (Theorem 4.4 of [53]) that gives a lower
bound on the Heegaard genus of a Haken manifold in terms of
the genus of a properly embedded incompressible surface and
the Heegaard genus of its complementary pieces. (See also
Proposition 23.40 of [29] for a similar result.) The maximal
number of properly embedded disjoint non-parallel essential
annuli in the complementary pieces also appears in the formula.
We now investigate this.

We claim that in a compact orientable atoroidal irreducible 3-manifold $N$,
there can be no more that $3|\chi(\partial N)|$ properly
embedded disjoint non-parallel essential annuli. Pick a
maximal such collection $A$. Suppose that two components of
$\partial A$ are parallel in $\partial N$. We may find two
that have no component of $\partial A$ in the annulus $A'$
between them. If these lie in distinct components of $A$,
say $A_1$ and $A_2$, then $A_1 \cup A' \cup A_2$ forms
an annulus which must be parallel to a component $A_3$ of
$A$. Then $A_1 \cup A_2 \cup A_3$ separates off a solid
torus in $N$. This
has a natural Seifert fibred structure in which $A_1$,
$A_2$ and $A_3$ are a union of fibres. Similarly, if the
components of $\partial A'$ lie in a single component $A_1$
of $A$, then $A_1$ separates off a Seifert fibred solid
torus in $N$ with one exceptional fibre,
since $N$ is irreducible and atoroidal. The union of
all such Seifert fibred solid tori is a Seifert fibre space $V$.
Each component must have base space that is a disc, and have
at most one exceptional fibre, since $N$ is atoroidal.
A simple counting argument then gives that at least half the
components of $A \cap V$ lie in $\partial V$. Let $A_-$
denote these annuli. Then $\partial A_-$ forms a collection
of curves in $\partial N$ with the property that no three
curves are parallel. There can be at most 
$3 |\chi(\partial N)|$ such curves. So there can be at most
${3 \over 2} |\chi(\partial N)|$ components of $A_-$ and so at most
$3 |\chi(\partial N)|$ components of $A$, which proves the
claim.

Theorem 4.4 of [53] gives a lower bound on $g(M_i, \partial M_i)$,
which is the minimal genus of a Heegaard surface for $M_i$,
subject to the condition that all boundary components for
$M_i$ lie on one side of this surface. We wish to relate this
to $\chi_-^h(M_i)$. Now, we may assume
that $S$ intersects each component of $\partial M$, for otherwise,
$M_i$ contains at least $i$ boundary components, and so
$\chi_-^h(M_i) \geq i -2$, which implies that $\{M_i \rightarrow M \}$
has positive Heegaard gradient. But when $S$ intersects each component 
of $\partial M$, there is a uniform upper bound ($t$, say) on
the number of components of $\partial M_i$. From any Heegaard surface for $M_i$,
we may therefore create one that has all components of $\partial M_i$
on one side, increasing the genus of the surface by at most
$t$. So, denoting the Heegaard genus of $M_i$ by
$g(M_i)$, we deduce
$$\chi_-^h(M_i) = 2 g(M_i) - 2 \geq 2g(M_i, \partial M_i) - 2 - 2t.$$

Now, for any $i \geq n$, consider the manifold obtained from $M_i$ by cutting
along $S_n \cup S_{2n} \cup \dots \cup S_{n \lfloor {i \over n} \rfloor}$.
We obtain $(\lfloor {i \over n} \rfloor - 1)$ copies of
$M(n)$ and one copy of $M(i - \lfloor {i \over n} \rfloor n + n)$.
Applying Theorem 4.4 of [53], we obtain
$$\eqalign{
\chi^h_-(M_i) &\geq 2 g(M_i, \partial M_i) - 2  - 2t\cr
& \geq {1 \over 7} \chi_-^h(M_i - (S_n \cup S_{2n} \cup
\dots S_{n \lfloor {i \over n} \rfloor})) + 2\cr
& \quad - 4 |\chi(S_n \cup S_{2n} \cup \dots \cup S_{n \lfloor {i \over n}
\rfloor})| - 2 - 2t \cr
& \geq \left \lfloor {i \over n} \right \rfloor |\chi(S)| - 2t
\cr
& \geq i {|\chi(S)| \over 2 n },}$$
when $i$ is sufficiently large.
So the Heegaard gradient of $\{ M_i \rightarrow M \}$ is
at least $|\chi(S)|/ 2 n$, which is non-zero. 
$\square$

\vskip 18pt
\centerline {\caps References}
\vskip 6pt

\item{1.} {\caps N. Alon}, {\sl Eigenvalues and expanders},
Combinatorica {\bf 6} (1986) 83--96.
\item{2.} {\caps N. Alon and V. D. Milman}, {\sl $\lambda_1$,
isoperimetric inequalities for graphs and superconcentrators},
J. Comb. Th. {\bf B 38} (1985) 78--88.
\item{3.} {\caps M. Berger and B. Gostiaux}, {\sl Differential
Geometry: Manifolds, curves and surfaces}, GTM 115, Springer-Verlag
(1988).
\item{4.} {\caps R. Brooks}, {\sl The spectral geometry of a 
tower of coverings}, J. Differential Geom. {\bf 23} (1986) 97--107.
\item{5.} {\caps P. Buser}, {\sl A note on the isoperimetric
constant}, Ann. Sci. Ecole Norm. Sup. {\bf 15} (1982) 213--230.
\item{6.} {\caps A. Casson and C. Gordon}, {\sl Reducing Heegaard splittings},
Topology and its Appl. {\bf 27} (1987) 275--283.
\item{7.} {\caps A. Casson and D. Jungreis}, {\sl Convergence groups 
and Seifert fibered $3$-manifolds,} Invent. Math.
{\bf 118} (1994) 441--456.
\item{8.} {\caps J. Cheeger}, {\sl A lower bound for the smallest
eigenvalue of the Laplacian}, Problems in Analysis (Ganning, ed.)
Princeton Univ. Press (1970) 195--199.
\item{9.} {\caps D. Cooper and D. Long}, 
{\sl Virtually Haken Dehn-filling}, J. Differential Geom. {\bf 52}
(1999) 173--187. 
\item{10.} {\caps J. Dodziuk}, {\sl Difference equations,
isoperimetric inequalities and transience of certain random
walks}, Trans. A.M.S. {\bf 284} (1984) 787--794.
\item{11.} {\caps J. Douglas}, {\sl The mapping theorem
of Koebe and the problem of Plateau}, J. Math. Phys. {\bf 10}
(1930-31) 106--130.
\item{12.} {\caps J. Douglas}, {\sl Solution of the
problem of Plateau}, Trans. A.M.S. {\bf 33} (1931) 263--321.
\item{13.} {\caps N. Dunfield and W. Thurston}, {\sl
The virtual Haken conjecture: experiments and examples}, 
Geom. Topol. {\bf 7} (2003) 399--441.
\item{14.} {\caps D. Epstein and R. Penner}, 
{\sl Euclidean decomposition
of non-compact hyperbolic manifolds}, J. Differential Geom. {\bf 27} (1988) 
67-80.
\item{15.} {\caps W. Floyd and A. Hatcher},  
{\sl Incompressible surfaces in punctured-torus bundles},
Topology Appl. {\bf 13} (1982) 263--282. 
\item{16.} {\caps M. Freedman, J. Hass and P. Scott,}
{\sl Least area incompressible surfaces in $3$-manifolds.}
Invent. Math. {\bf 71} (1983) 609--642.
\item{17.} {\caps C. Frohman}, {\sl The topological uniqueness
of triply periodic minimal surfaces in ${\Bbb R}^3$}, J. Differential
Geom. {\bf 31} (1990) 271--283.
\item{18.} {\caps D. Gabai}, {\sl Convergence groups are Fuchsian groups},
Ann. of Math. (2) {\bf 136} (1992) 447--510. 
\item{19.} {\caps D. Gabai}, {\sl Homotopy hyperbolic $3$-manifolds 
are virtually hyperbolic}, J. Amer. Math. Soc. {\bf 7}
(1994) 193--198.
\item{20.} {\caps D. Gabai}, {\sl On the geometric and topological rigidity of 
hyperbolic \break $3$-manifolds}, J. Amer. Math. Soc. {\bf 10} (1997) 37--74.
\item{21.} {\caps W. Haken}, {\sl Some results on surfaces in
3-manifolds}, Studies in Modern Topology (1968) 34--98.
\item{22.} {\caps R. Hamilton,} {\sl 
Three-manifolds with positive Ricci curvature,}
J. Differential Geom. {\bf 17} (1982) 255--306.
\item{23.} {\caps R. Hamilton,} {\sl
The formation of singularities in the Ricci flow.}
Surveys in differential geometry, Vol. II (Cambridge, MA, 1993), 
7--136.
\item{24.} {\caps R. Hamilton,} {\sl The Harnack estimate for the Ricci flow.}
J. Differential Geom. {\bf 37} (1993) 225--243.
\item{25.} {\caps R. Hamilton,} {\sl 
Non-singular solutions of the Ricci flow on three-manifolds,}
Comm. Anal. Geom. {\bf 7} (1999) 695--729.
\item{26.} {\caps J. Hass, S. Wang and Q. Zhou}, {\sl On
finiteness of the number of boundary slopes of immersed
surfaces in 3-manifolds}, Proc. Amer. Math. Soc. {\bf 130} (2002) 1851--1857. 
\item{27.} {\caps J. Hempel}, {\sl 3-Manifolds}, Ann. of Math. Studies, No. 86,
Princeton Univ. Press, Princeton, N. J. (1976)
\item{28.} {\caps K. Ichihara}, {\sl Heegaard gradient of Seifert fibered
3-manifolds},  Bull. London Math. Soc.  {\bf 36}  (2004)  537--546. 
\item{29.} {\caps K. Johannson}, {\sl Topology and combinatorics of
3-manifolds}, Springer-Verlag (1995).
\item{30.} {\caps D. Kazhdan}, {\sl Connection of the dual space of a
group with structure of its closed subgroups}, Functional Anal.
Appl. {\bf 1} (1967) 63-65.
\item{31.} {\caps M. Lackenby}, {\sl Word hyperbolic Dehn surgery},
Invent. Math. {\bf 140} (2000) 243--282.
\item{32.} {\caps M. Lackenby}, {\sl Classification of
alternating knots with tunnel number one,} Comm. Anal.
Geom. {\bf 13} (2005) 151-186 .
\item{33.} {\caps C. Leininger}, {\sl 
Surgeries on one component of the Whitehead link are virtually
fibered}, Topology {\bf 41} (2002) 307--320. 
\item{34.} {\caps G. Livesay}, {\sl Fixed point free involutions
on the 3-sphere}, Ann. Math {\bf 72} (1960) 603--611.
\item{35.} {\caps A. Lubotzky}, {\sl Discrete Groups,
Expanding Graphs and Invariant Measures}, Progress
in Math. {\bf 125} (1994)
\item{36.} {\caps A. Lubotzky}, {\sl Eigenvalues of the 
Laplacian, the first Betti number and the congruence subgroup
problem,} Ann. Math. {\bf 144} (1996) 441--452.
\item{37.} {\caps A. Lubotzky and R. Zimmer}, {\sl
Variants of Kazhdan's property for subgroups of semisimple groups}, 
Israel J. Math. {\bf 66} (1989) 289--299.
\item{38.} {\caps W. Meeks III and S. Yau},
{\sl Topology of three-dimensional manifolds and 
the embedding problems in minimal surface theory},
Ann. Math. {\bf 112} (1980) 441--484. 
\item{39.} {\caps W. Meeks III, L. Simon and S. Yau,}
{\sl Embedded minimal surfaces, exotic spheres, and manifolds 
with positive Ricci curvature,}
Ann. Math. {\bf 116} (1982) 621--659.
\item{40.} {\caps J. Morgan}, {\sl The Smith conjecture.} 
Pure Appl. Math. {\bf 112} (1984).
\item{41.} {\caps Y. Moriah}, {\sl On boundary primitive manifolds
and a theorem of Casson-Gordon}, Topology Appl. {\bf 125} (2002) 571--579.
\item{42.} {\caps G. Perelman,} {\sl The entropy formula for the 
Ricci flow and its geometric applications,} Preprint,
available at \hfill\break
http://front.math.ucdavis.edu/math.DG/0211159
\item{43.} {\caps G. Perelman,} {\sl  Ricci flow with surgery on three-manifolds,}
Preprint, available at http://front.math.ucdavis.edu/math.DG/0303109
\item{44.} {\caps G. Perelman,} {\sl Finite extinction time for the 
solutions to the Ricci flow on certain three-manifolds,} Preprint,
available at \hfill\break
http://front.math.ucdavis.edu/math.DG/0307245 
\item{45.} {\caps J. Pitts and J. H. Rubinstein,} {\sl Existence 
of minimal surfaces of bounded topological type in three-manifolds.}
Miniconference on geometry and partial differential equations 
(Canberra, 1985), 163--176.
\item{46.} {\caps T. Rad\'o}, {\sl The problem of least
area and the problem of Plateau}, Math. Z. {\bf 32} (1930)
763--796.
\item{47.} {\caps T. Rad\'o}, {\sl On Plateau's problem},
Ann. Math. {\bf 31} (1930) 457--469.
\item{48.} {\caps A. Reid}, {\sl A non-Haken hyperbolic 3-manifold 
covered by a surface bundle}, Pacific J. Math. {\bf 167} (1995) 163--182.
\item{49.} {\caps J. H. Rubinstein}, {\sl Polyhedral minimal surfaces,
Heegaard splittings and decision problems for 3-dimensional
manifolds}, Proceedings of the Georgia Topology Conference,
AMS/IP Stud. Adv. Math, vol. 21, Amer. Math Soc. (1997) 1--20.
\item{50.} {\caps J. H. Rubinstein}, {\sl Minimal surfaces
in geometric 3-manifolds}, Preprint.
\item{51.} {\caps M. Scharlemann and A. Thompson}, {\sl
Thin position for $3$-manifolds,} Geometric topology (Haifa,
1992), 231--238, Contemp. Math., 164.
\item{52.} {\caps R. Schoen and S. Yau}, {\sl Existence of
incompressible surfaces and the topology of 3-manifolds
with non-negative scalar curvature}, Ann. Math. {\bf 119}
(1979) 127--142.
\item{53.} {\caps J. Schultens}, {\sl Heegaard genus formula for
Haken manifolds}, Preprint, available at http://front.math.ucdavis.edu/math.GT/0108028
\item{54.} {\caps P. Scott}, {\sl There are no fake Seifert fibre spaces 
with infinite $\pi \sb{1}$}, Ann. of Math. (2) {\bf 117} (1983) 35--70. 
\item{55.} {\caps M. Stocking}, {\sl Almost normal surfaces in 3-manifolds},
Trans. Amer. Math. Soc. {\bf 352} (2000) 171--207.
\item{56.} {\caps R. Tanner}, {\sl Explicit concentrators from
generalised $N$-gons}, SIAM J. Alg. Discr. Meth. {\bf 5}
(1984) 287--294.
\item{57.} {\caps W. Thurston}, {\sl The geometry and topology
of 3-manifolds (Lecture Notes)}, Princeton (1980).
\item{58.} {\caps W. Thurston}, {\sl Three-dimensional manifolds, 
Kleinian groups and
hyperbolic geometry}, Bull. Amer. Math. Soc. {\bf 6} (1982) 357-381.
\item{59.} {\caps B. White}, {\sl The space of minimal submanifolds
for varying Riemannian metrics}, Indiana Math. Journal {\bf 40} 
(1991) 161--200.

\vskip 12pt

\+ Mathematical Institute, Oxford University, \cr
\+ 24-29 St Giles', Oxford OX1 3LB, United Kingdom. \cr

\end